\documentclass[12pt]{article}

\usepackage{amsmath}
\usepackage{amsthm,amscd,amsfonts}
\usepackage{amssymb, upref, color}
\usepackage{amsmath,amssymb,amsthm}
\usepackage{color}
\usepackage[colorlinks]{hyperref}
\usepackage{graphicx}
\usepackage{epsf,epsfig,subfigure, verbatim}
\usepackage{latexsym,bm}

\numberwithin{equation}{section}

\theoremstyle{plain}
\newtheorem{exam}{Example}[section]
\newtheorem{theorem}[exam]{Theorem}
\newtheorem{lemma}[exam]{Lemma}
\newtheorem{remark}[exam]{Remark}

\newtheorem{proposition}[exam]{Proposition}
\newtheorem{definition}[exam]{Definition}
\newtheorem{corollary}[exam]{Corollary}

\begin{document}
\date{}


\title{ Linear Quaternion Differential Equations: Basic Theory and Fundamental Results
}
\author{ Kit Ian Kou$^{1,2}$\footnote{ Kit Ian Kou acknowledges financial support from the National Natural Science Foundation of China under Grant (No. 11401606), University of Macau (No. MYRG2015-00058-FST and No. MYRG099(Y1-L2)-FST13-KKI) and the Macao Science and Technology
Development Fund (No. FDCT/094/2011/A and No. FDCT/099/2012/A3). }
 \,\,\,\,\,  Yong-Hui Xia$^{1,3}$
 \footnote{ Corresponding author. Email: xiadoc@163.com.  Yonghui Xia was supported by the National Natural
Science Foundation of China under Grant (No. 11671176 and No. 11271333), Natural
Science Foundation of Zhejiang Province under Grant (No. Y15A010022), Marie Curie Individual Fellowship within the European Community Framework Programme(MSCA-IF-2014-EF, ILDS - DLV-655209), the Scientific Research Funds of Huaqiao University and China Postdoctoral Science Foundation (No. 2014M562320). }
\\
{\small 1. School of Mathematical Sciences, Huaqiao University, 362021, Quanzhou, Fujian, China.}\\
{\small\em xiadoc@163.com; xiadoc@hqu.edu.cn (Y.H.Xia)}\\
{\small 2.Department of Mathematics, Faculty of Science and Technology, University of Macau, Macau}\\
{\small \em  kikou@umac.mo (K. I. Kou)  }
 \\
 {\small 3. Department of Mathematics, Zhejiang Normal University, Jinhua, 321004, China}\\
}
 \maketitle


\begin{center}
\begin{minipage}{140mm}
\begin{abstract}

Quaternion-valued differential equations (QDEs) is a new kind of differential equations which have many applications in physics and life sciences. The largest difference between QDEs and ODEs is the algebraic structure. Due to the non-commutativity of the quaternion algebra, the algebraic structure of the solutions to the QDEs is completely different from ODEs. {\em It is actually a right free module, not a linear vector space.}

This paper establishes {\bf a systematic frame work for the theory of linear QDEs}, which can be applied to quantum mechanics, fluid mechanics, Frenet frame in differential geometry, kinematic modelling, attitude dynamics, Kalman filter design and spatial rigid body dynamics, etc. We prove that the algebraic structure of the solutions to the QDEs is actually a right free module, not a linear vector space. On the non-commutativity of the quaternion algebra, many concepts and properties for the ordinary differential equations (ODEs) can not be used. They should be redefined accordingly.
A definition of {\em Wronskian} is introduced under the framework of quaternions which is different from standard one in the ordinary differential equations. Liouville formula for QDEs is given. Also, it is necessary to treat the eigenvalue problems with left- and right-sides, accordingly. Upon these, we studied the solutions to the linear QDEs. Furthermore, we present two algorithms to evaluate the fundamental matrix. Some concrete examples are given to show the feasibility of the obtained algorithms. Finally, a conclusion and discussion ends the paper.

\end{abstract}



\end{minipage}
\end{center}
\section{\bf Introduction and Motivation}


The ordinary differential equations (ODEs) have been well studied. Theory of ODEs is very systematic and rather complete (see monographs, e.g. \cite{Arnold, Chicone,Chow,DingLi,Hale,Hartman,ZWN}). Quaternion-valued differential equations (QDEs) is a new kind of differential equations which have many applications in physics and life sciences. Recently, there are few papers trying to study QDEs \cite{Mawhin,Zhang1,Leo3,W,Zhang2}. Up till now, the theory of the QDEs is not systemic. Many questions remains open. For example, even the structure of the solutions to linear QDEs is unknown. What kind of vector space is the solutions to linear QDEs? What is the difference between QDEs and ODEs? First of all, why should we study the QDEs? In the following, we will introduce the background of the quaternion-valued differential equations?

\subsection{Background for QDEs}

Quaternions are 4-vectors whose multiplication rules are governed by a simple non-commutative division algebra. We denote the quaternion $q=(q_0, q_1, q_2, q_3)^T\in \mathbb{R}^4$ by
\[
q= q_0 + q_1 \mathbf{i} + q_2 \mathbf{j} + q_3 \mathbf{k},
\]
where $q_0, q_1, q_2, q_3$ are real numbers and $\mathbf{i},\mathbf{j},\mathbf{k}$ satisfy the multiplication table formed by
\[
\mathbf{i}^2=\mathbf{j}^2=\mathbf{k}^2=\mathbf{ijk}=-1,\,\,  \mathbf{ij}=-\mathbf{ji}=\mathbf{k},\,\, \mathbf{ki}=-\mathbf{ik}=\mathbf{j},\,\,\,\mathbf{jk}=-\mathbf{kj}=\mathbf{i}.
\]
The concept was originally invented by Hamilton in 1843 that extends the complex numbers to four-dimensional space.

\begin{figure}[h]
\centering
\includegraphics[width=3.8in,height=2.8in]{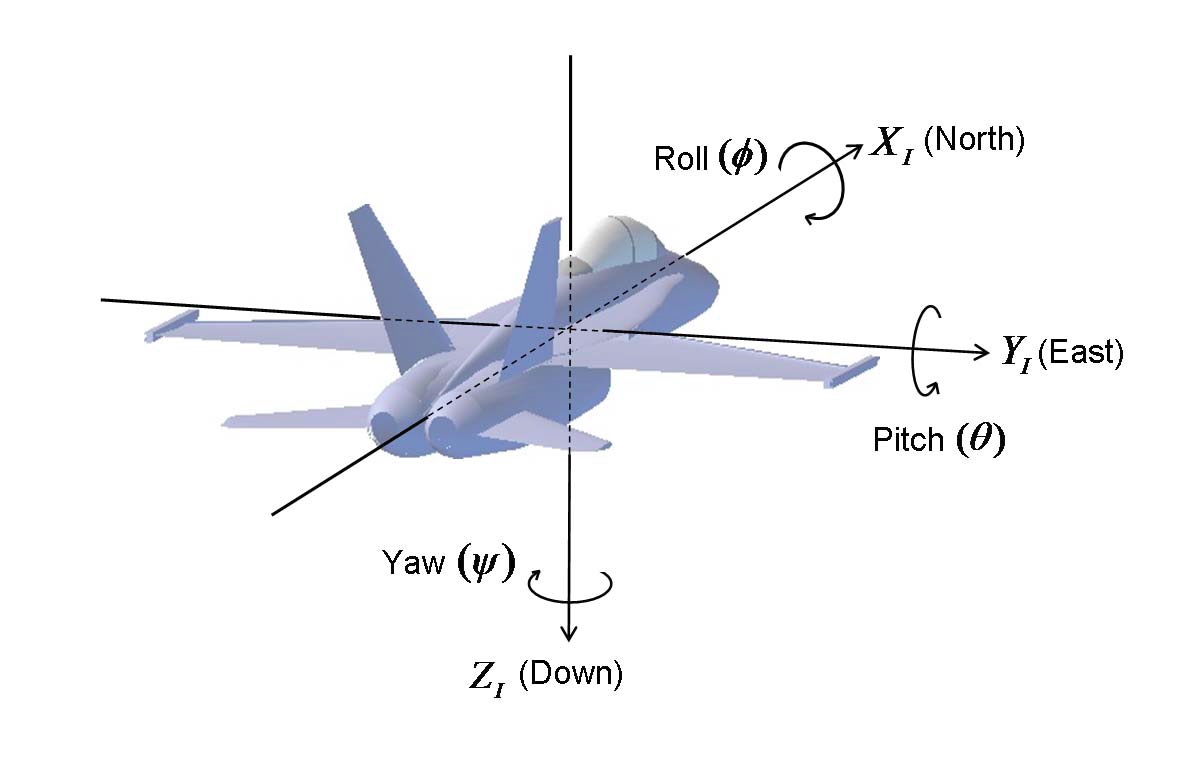}
\caption{Inertial Frame} \label{fig1}
\end{figure}

\begin{figure}[h]
\centering
\includegraphics[width=3.8in,height=2.2in]{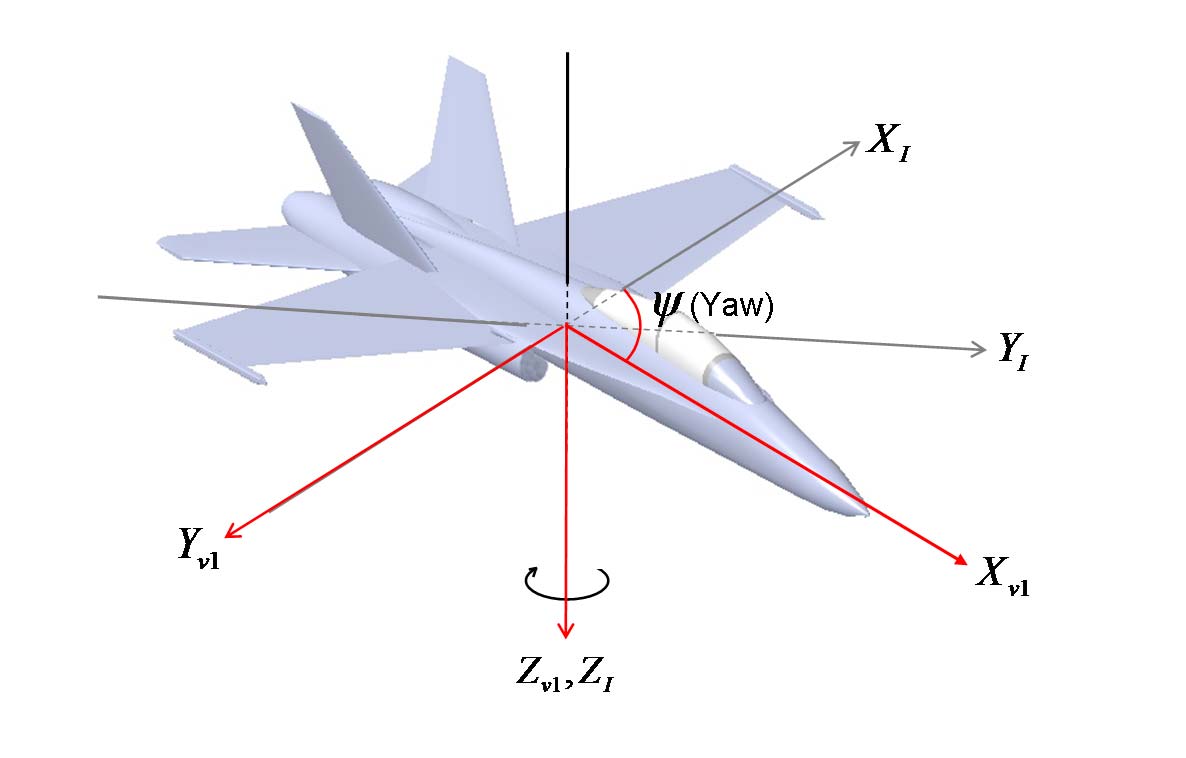}
\caption{  Vehicle-1 Frame.} \label{fig2}
\end{figure}

\begin{figure}[h]
\centering
\includegraphics[width=3.8in,height=2.2in]{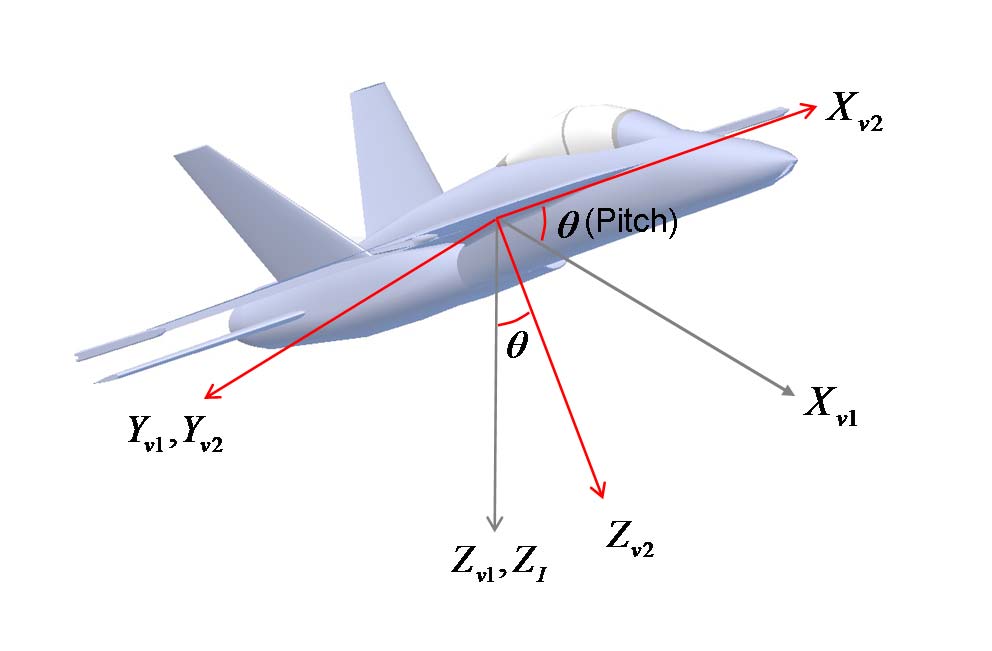}
\caption{  Vehicle-2 Frame.} \label{fig3}
\end{figure}

\begin{figure}[h]
\centering
\includegraphics[width=3.8in,height=2.2in]{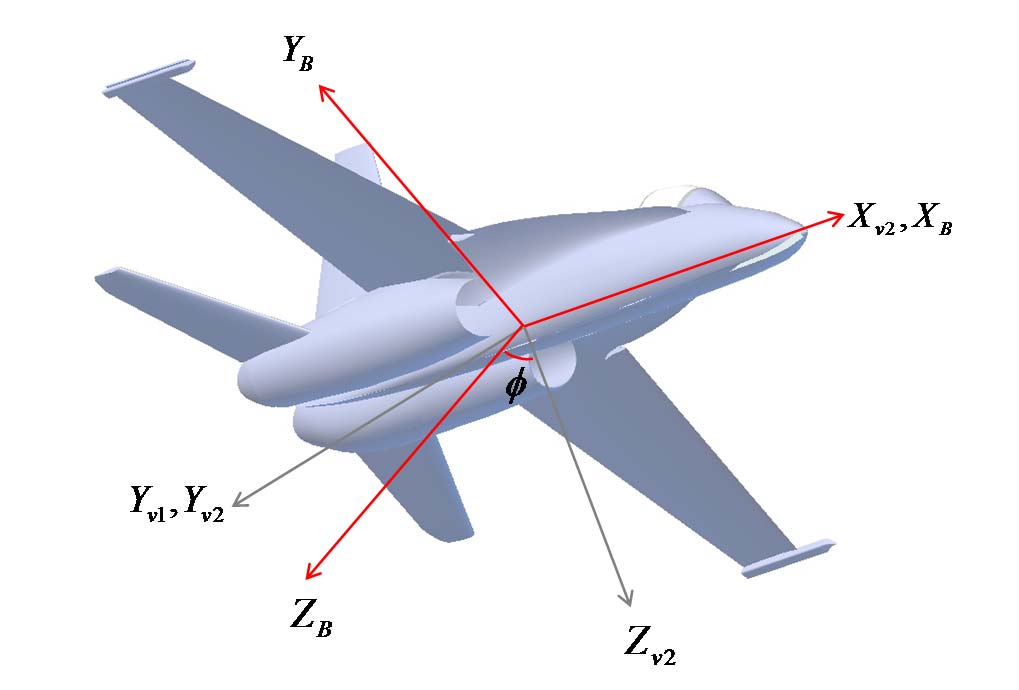}
\caption{  Body Frame.} \label{fig4}
\end{figure}

Quaternions have shown advantages over real-valued vectors in physics and engineering applications for their powerful modeling of rotation and orientation. Orientation can be defined as a set of parameters
that relates the angular position of a frame to another
reference frame. There are numerous methods for
describing this relationship. Some are easier to visualize
than the others. Each has some kind of limitations.
Among them, 
 Euler angles and
quaternions are commonly used. 
We give an example from \cite{AHS, AHS1} for illustration. Attitude and Heading Sensors (AHS) from CH robotics can provide orientation information using both Euler angles and quaternions. Compared to quaternions, Euler angles are simple and intuitive
  (see Figure \ref{fig1},\ref{fig2},\ref{fig3} and \ref{fig4}). The complete rotation matrix for moving from the inertial frame to the body frame is given by the multiplication of three matrix taking the form
\[
R_{I}^{B}(\phi,\theta,\psi)=R_{V_2}^{B}(\phi)R_{V_1}^{V_2}(\theta)R_{I}^{V_1}(\psi),
\]
where $R_{V_2}^{B}(\phi),R_{V_1}^{V_2}(\theta),R_{I}^{V_1}(\psi)$ are the rotation matrix moving from vehicle-2 frame to the body frame,
 vehicle-1 frame to vehicle-2 frame, inertial frame to vehicle-1 frame, respectively.
\begin{figure}[h]
\centering
\includegraphics[width=3.2cm,height=2.2in]{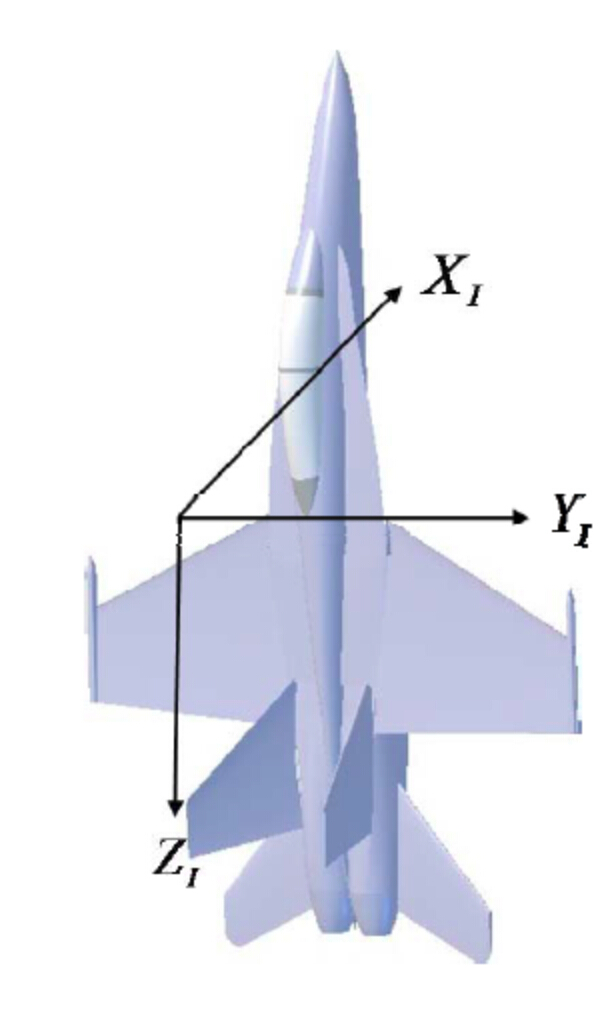}
\caption{Gimbal Lock.} \label{fig5}
\end{figure}

On the other hand, Euler angles are limited by a phenomenon called ``Gimbal Lock". The cause of such phenomenon is that when the pitch angle is 90 degrees
(see Figure \ref{fig5}). An orientation sensor
 that uses Euler angles will always fail to produce reliable estimates when the pitch angle approaches 90 degrees. 
  This is a serious shortcoming of Euler angles and can only be solved by switching to a different representation method. Quaternions provide an alternative measurement technique that does not suffer from gimbal lock. Therefore, all CH Robotics attitude sensors use quaternions so that the output is always valid even when Euler Angles are not.

The attitude quaternion estimated by CH Robotics orientation sensors encodes rotation from the ``inertial frame" to the sensor ``body frame."  The inertial frame is an Earth-fixed coordinate frame defined so that the x-axis points north, the y-axis points east, and the z-axis points down as shown in Figure 1.  The sensor body-frame is a coordinate frame that remains aligned with the sensor at all times.  Unlike Euler angle estimation, only the body frame and the inertial frame are needed when quaternions are used for estimation.

Let the vector $q=(q_0, q_1, q_2, q_3)^{T}$ be defined as the unit-vector quaternion encoding rotation from the inertial frame to the body frame of the sensor,
where $T$ is the vector transpose operator. The elements $q_1, q_2$, and $q_3$ are the ``vector part" of the quaternion, and can be thought of as a vector about which rotation should be performed. The element $q_0$ is the ``scalar part" that specifies the amount of rotation that should be performed about the vector part.  Specifically, if $\theta$ is the angle of rotation and the vector $e=(e_x,e_y,e_z)$ is a unit vector representing the axis of rotation, then the quaternion elements are defined as
\[
q
=
\left ( \begin{array}{c}
q_0
\\
q_1
\\
q_2
\\
q_3
\end{array}
\right )
=\left ( \begin{array}{c}
\cos (\theta/2) 
\\
e_x \sin (\theta/2) 
\\
e_y \sin (\theta/2) 
\\
e_z \sin (\theta/2) 
\end{array}
\right ).
\]

In fact, Euler's
theorem states that given two coordinate systems, there is one invariant axis (namely, Euler axis), along which
measurements are the same in both coordinate systems. Euler's
theorem also shows that it is possible to move
from one coordinate system to the other through one rotation $\theta$ about that invariant axis. Quaternionic representation of the attitude is based on Euler's
theorem. Given a unit vector
$e=(e_x,e_y,e_z)$ along the Euler axis, the quaternion is defined to be
\[
q=\left ( \begin{array}{c}
\cos (\theta/2) 
\\
e \sin (\theta/2) 
\end{array}
\right ).
\]

\begin{figure}[h]
\centering
\includegraphics[width=2in,height=2.5in]{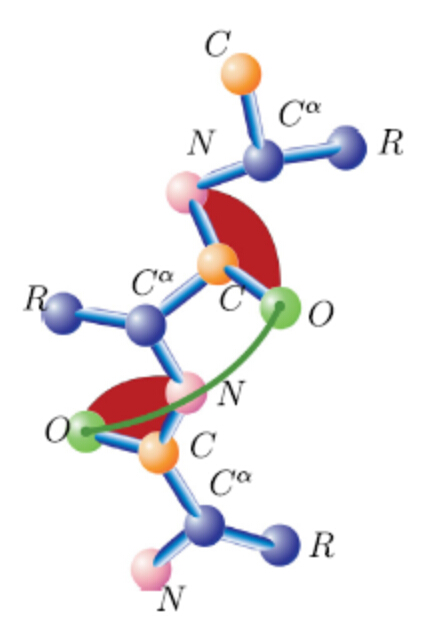}
\caption{The transformation that sends one peptide plane to the next is a screw
motion.} \label{fig6}
\end{figure}

Besides the attitude orientation, quaternions have been widely applied to study life science, physics and engineering. For instances,
 an interesting application is the description of protein structure (see Figure \ref{fig6}) \cite{Protein}, neural networks \cite{NN}, spatial rigid body transformation \cite{Body}, Frenet frames in differential geometry, fluid mechanics, quantum mechanics, and so on.

\subsubsection {  Quaternion Frenet frames in differential geometry}

\begin{figure}[h]
\centering
\includegraphics[width=2.8in,height=1.8in]{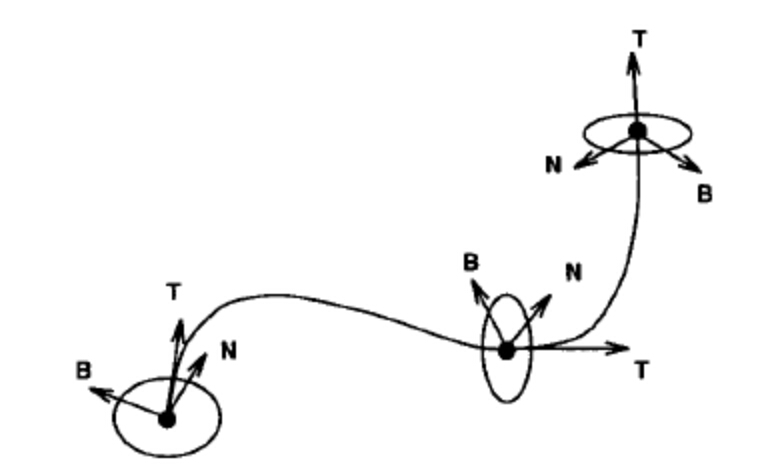}
\caption{Frenet frame for a curve.} \label{fig7}
\end{figure}

The differential geometry of curves \cite{DG1,DG2,DG3} traditionally begins with a vector $\vec{x}(s)$ that describes the curve parametrically as a function
of $s$ that is at least three times differentiable. Then the tangent vector
$\vec{T}(s)$ is well-defined
at every point $\vec{x}(s)$ and one may choose two additional orthogonal vectors in the plane
perpendicular to $\vec{T}(s)$ to form a complete local orientation frame. Provided the curvature of $\vec{x}(s)$ vanishes nowhere, we can choose this local coordinate system to be the
Frenet frame consisting of the tangent $\vec{T}(s)$,
the binormal $\vec{B}(s)$, and the principal normal $\vec{N}(s)$ (see Figure \ref{fig7}), which are given in terms of the
curve itself by these expressions:
\[
\begin{array}{lll}
\vec{T}(s)&=& \displaystyle\frac{\vec{x}'(s)}{\|\vec{x}'(s)\|},
\\
\vec{B}(s)&=& \displaystyle\frac{\vec{x}'(s)\times \vec{x}''(s)}{\|\vec{x}'(s)\times \vec{x}''(s)\|},
\\
\vec{N}(s)&=&  \displaystyle\vec{B}(s) \times \vec{T}(s),
\end{array}
\]
The standard frame configuration is illustrated in Figure 7.
 The Frenet frames obeys the following differential equations.
\begin{equation}
\left ( \begin{array}{cccc}
\vec{T}'(s)
\\
\vec{B}'(s)
  \\
\vec{N}'(s)
\end{array}
\right )
=
v(s)\left ( \begin{array}{cccc}
0 & \kappa(s) & 0
\\
- \kappa(s) & 0 &  \tau(s)
  \\
0  & -\tau(s) & 0
\end{array}
\right )
\left ( \begin{array}{cccc}
\vec{T}(s)
\\
\vec{B}(s)
  \\
\vec{N}(s)
\end{array}
\right )
\label{F1}
\end{equation}
Here $v(s) =\|\vec{x}'(s)\|$ is the scalar magnitude of the curve derivative, and the intrinsic
geometry of the curve is embodied in the curvature $\kappa(s)$ and the torsion $\tau(s)$, which
may be written in terms of the curve itself as
\[
\begin{array}{lll}
\kappa(s)(s)&=& \displaystyle\frac{\|\vec{x}'(s)\times \vec{x}''(s)\|}{\|\vec{x}'(s)\|^3},
\\
\tau(s)(s)&=& \displaystyle\frac{(\vec{x}'(s)\times \vec{x}''(s))\cdot \vec{x}'''(s)}{\|\vec{x}'(s)\times \vec{x}''(s)\|^2}.
\end{array}
\]
Handson and Ma \cite{Handson} pointed out that all $3D$ coordinate frame can be express in the form of quaternion.
Then $q'(t)=(q'_0(t),q'_1(t),q'_2(t),q'_3(t))$ can be written as a first-order quaternion differential equation
\[
q'(t)=\frac{1}{2}vK  q(t)
\]
with
\[
K=\left ( \begin{array}{cccc}
0 & -\tau & 0 & -\kappa
\\
\tau &0 &  \kappa &0
  \\
0  & -\kappa & 0 & \tau
    \\
\kappa & 0 & -\tau & 0
\end{array}
\right ).
\]

\subsubsection {  QDEs appears in kinematic modelling and attitude dynamics}
Global Positioning System (GPS) have developed rapidly in recent years. It is compulsory for the navigation and guidance.
In fact, GPS is based on attitude determination. And attitude determination
results in the measurement of vehicle orientation with respect to a reference frame. 
 Wertz \cite{Wertz} derived the general attitude motion equations for the
quaternions assuming a constant rotation rate over an infinitesimal time $\Delta t$ such that
\[
q_{k+1}=q_{k}+\dot{q}\Delta t
\]
where $q =(q_1,q_2,q_3,q_4)^T$, and the quaternion differentiated with respect to time is
approximated by
\[
\dot {q}=\frac{1}{2}\Omega \cdot q_k
\]
with the skew symmetric form of the body rotations about the reference frame as
\[
\Omega=\left ( \begin{array}{cccc}
0 & \omega_z & -\omega_y & \omega_x
\\
- \omega_z &0 &  \omega_x & \omega_y
  \\
\omega_y   & -\omega_x & 0 & \omega_z
    \\
-\omega_x & -\omega_y & -\omega_z & 0
\end{array}
\right )
\]
and $ \omega=(\omega_x,  \omega_y, \omega_z)^T$ is the angular velocity of body rotation.

Another attitude dynamic is the orientation tracking for humans and robots
using inertial sensors. In \cite{Marins}, the authors
designed a Kalman filter for estimating
orientation. They developed a process model
of a rigid body under rotational motion.  The state vector of a process model will consists of the
angular rate $\omega$ and parameters for characterizing orientation $n$. For human body
motions, it is shown that quaternion rates are related to angular rates
through the following first order quaternion differential equation:
\[
\dot {n}=\frac{1}{2}n \, \omega,
\]
where $\omega$ is treated as a quaternion with
zero scalar part.

Moreover, a general formulation for rigid body rotational
dynamics is developed using quaternions. 
Udwadia et al. \cite{Udwadia} provided a simple and direct route for obtaining
Lagrange's equation describing rigid body rotational motion in
terms of quaternions. He presented a second order quaternionic differntial equation as follows.
\[
4E^TJE\ddot{u} + 8 \dot{E}^TJE\dot{u} +4 J_0N(\dot{u})u = 2E^T \Gamma_B- 2(\omega^T J \omega)u,
\]
where $u=(u_0,u_1,u_2,u_3)^T$ is the unit quaternion with $u^Tu=1$ and $E$ is the orthogonal matrix given by
\[
E=\left ( \begin{array}{cccc}
u_0 & u_1 & u_2 & u_3
\\
- u_1 &u_0 &  u_3 & -u_2
  \\
-u_2   & -u_3 & u_0 & u_1
    \\
-u_3 & u_2 & -u_1 & u_0
\end{array}
\right ),
\]
$\Gamma_B$ is the four-vector containing the components of the
body torque about the body-fixed axes. $\omega=(0,\omega_1,\omega_2,\omega_3)^T$, $\omega_1,\omega_2,\omega_3$ are components of angular velocity.
The inertia matrix of the rigid body is represented by $\hat{J} =diag(J_1,J_2,J_3)^T$
and the $4\times 4$ diagonal matrix $J=diag(J_0,J_1,
J_2,J_3)^T$, where $J_0$ is any positive number.

\subsubsection {  QDE appears in fluid mechanics}
Recently Gibbon et
al. \cite{Gibbon1,Gibbon2} used quaternion algebra to study the fluid mechanics.
They proved that the three-dimensional incompressible Euler equations have a natural {\em quaternionic Riccati
structure} in the dependent variable. He proved that the 4-vector $\zeta$ for the three-dimensional incompressible Euler equations evolves according to the {\em quaternionic Riccati equation}
\[
\frac{D \zeta}{d t}+ \zeta^2 +\zeta_p=0,
\]
where $\zeta=(\alpha,\chi)^T$ and $\zeta_p=(\alpha_p,\chi_p)^T$ are a 4-vectors, $\alpha$ is a scalar and $\chi$ is a 3-vector, where $\alpha_p$, $\chi_p$ are defined in terms of the Hessian matrix $P$ \cite{Gibbon1}.
  To convince the reader that this structure is robust and no accident, it is shown that the more complicated equations for ideal magneto-hydrodynamics (MHD) can also be written in an
equivalent quaternionic form.
Roubtsov and Roulstone \cite{RR1,RR2} also expressed the
nearly geostrophic flow in terms of quaternions.

\subsubsection {  QDE appears in quantum mechanics}
Quaternion has been used to study the quantum mechanics for a very long time, and the quternionic quantum mechanic theory is well developed (see, for instance, Refs. \cite{Alder1,Alder2,Leo1,Leo2} and the references therein). In the standard formulation of nonrelativistic quantum mechanics, the complex wave function
$\Phi(x, t)$, describing a particle without spin subjected to the influence of a real potential $V(x, t)$,
satisfies the Schr$\ddot{o}$dinger equation
\[
  \partial_{t}\Phi(x, t)= \frac{\mathbf{i}}{{\hbar} } \Big[ \frac{ {\hbar}^2}{2m}\nabla^2 - V(x,t)\Big]\Phi(x, t).
\]
Alder generalized the classical mechanics to quaternionic quantum mechanics \cite{Alder2}, the anti-self-adjoint operator
\[
  {\cal{A}}^{V}(x, t)= \frac{\mathbf{i}}{{\hbar} } \Big[ \frac{ {\hbar}^2}{2m}\nabla^2 - V(x,t)\Big],
\]
can be generalized by introducing the complex potential $W(x, t) = |W(x, t)| \exp[\mathbf{i}\theta(x, t)]$,
\[
  {\cal{A}}^{V,W}(x, t)= \frac{\mathbf{i}}{{\hbar} } \Big[ \frac{ {\hbar}^2}{2m}\nabla^2 - V(x,t)\Big]+  \frac{\mathbf{j}}{{\hbar} } W(x, t) ,
\]
where the quaternionic imaginary units $\mathbf{i}, \mathbf{j}$ and $\mathbf{k}$ satisfy the following associative but noncommutative
algebra:
\[
\mathbf{i}^2=\mathbf{j}^2=\mathbf{k}^2=\mathbf{ijk}=-1.
\]
As pointed out in \cite{Leo1}, this generalization for the anti-self-adjoint Hamiltonian operator, the
quaternionic wave function $\Phi(x, t)$ satisfies the following equation:
\[
  \partial_{t}\Phi(x, t)= \Big\{\frac{\mathbf{i}}{{\hbar} } \Big[ \frac{ {\hbar}^2}{2m}\nabla^2 - V(x,t)\Big]+ \frac{\mathbf{j}}{{\hbar} } W(x, t)\Big\} \Phi(x, t).
\]
Later, \cite{Leo2} gave further study on the quaternionic quantum mechanics. The complex wave function
$\Phi(x, t)$, describing a particle without spin subjected to the influence of a potential $V(x, t)$,
satisfies the Schr$\ddot{o}$dinger equation
\[
\mathbf{i} \frac{ {\hbar}^2}{2m}\Phi_{xx}(x, t)-(\mathbf{i} V_1 + \mathbf{j} V_2 + \mathbf{k} V_3)\Phi(x, t)=\hbar \Phi_t(x, t).
\]
This partial differential equation can be reduced by using the well-known separation of variables,
\[
\Phi(x, t)=\phi(x)\exp\{-\mathbf{i} Et/{\hbar}\},
\]
to the following second order quaternion differential equation
\[
\mathbf{i} \frac{ {\hbar}^2}{2m}\ddot\phi(x)-(\mathbf{i} V_1 + \mathbf{j} V_2 + \mathbf{k} V_3)\phi(x)=- \phi(x)E \mathbf{i}.
\]


\subsection{History and motivation}

Though the above mentioned works show the existence of the quaternionic differential equation
structure appearing in differential geometry, fluid mechanics, attitude dynamics, quantum mechanics, etc. for the author knowledge, there is no research on the systematic study for the theory of the quaternion-valued differential equations. Upon the non-commutativity of the quaternion algebra, the study on the quaternion differential equations becomes a challenge and much involved. The studies are rare in this  area. However, if you would like to found more behavior of a complicated dynamical system, it is necessary to make further mathematical analysis.
By using the real matrix representation of left/right acting quaternionic
operators, Leo and Ducati \cite{Leo3}
tried to solve some simple second order quaternionic differential equations. Later, Campos and Mawhin \cite{Mawhin} studied the existence of periodic solutions of one-dimensional first order periodic {\em quaternion differential equation}. Wilczynski \cite{W} continued this study and payed more attention on the existence of two periodic solutions of quaternion Riccati equations. Gasull et al \cite{Zhang1} studied a one-dimensional quaternion autonomous homogeneous differential equation
\[
\dot{q} = aq^n.
\]
Their interesting results describe the phase portraits of the homogeneous quaternion differential equation.
For $n=2,3$, the existence of
 periodic orbits, homoclinic
loops, invariant tori and the first integrals are discussed in \cite{Zhang1}.
Zhang \cite{Zhang2} studied the global structure of the one-dimensional
quaternion Bernoulli equations
\[
\dot{q} = aq + aq^n.
\]
By using the Liouvillian theorem of integrability and the topological
characterization of 2-dimensional torus: orientable compact connected
surface of genus one, Zhang proved that the quaternion Bernoulli equations
may have invariant tori, which possesses a full Lebesgue measure subset
of $\mathbb{H}$. Moreover, if $n = 2$ all the invariant tori are full of periodic orbits;
if $n = 3$ there are infinitely many invariant tori fulfilling periodic orbits
and also infinitely many invariant ones fulfilling dense orbits.

As we know, the second order quaternionic differential equation can be easily transformed to a two-dimensional QDEs. In fact, for a second order quaternionic differential equation
\[
\frac{{\rm d^2} x}{{\rm} d t^2} + q_1(t) \frac{{\rm d}  x}{{\rm d} t } +q_2(t) x=0,
\]
by the transformation
\[
\left\{
\begin{array}{ccc}
x_1& =& x
\\
x_2 &= &\dot{x}
\end{array}
\right.
\]
then we have
\[
\left\{
\begin{array}{lll}
\dot{x}_1& =& x_1
\\
\dot{x}_2 &= & - q_2(t) x_1 - q_1(t) x_2 .
\end{array}
\right.
\]
Most of above mentioned works focused on the one-dimensional QDE and they mainly discussed the qualitative properties of the QDE (e.g. the existence of periodic orbits, homoclinic loops and invariant tori, integrability, the existence of periodic solutions and so on).
They did not provide the algorithm to compute the exact solutions to linear QDEs. The most basic problem is to solve the QDEs. Since quaternions are a non-commutative extension of the complex numbers, a lot of concepts and theory of ODEs are not valid for QDEs. Up till now, there is no systematic theory on the algorithm of solving QDEs. What is superposition principle? Can the solutions of QDEs span a linear vector space?  How is the structure of the general solutions? How to determine the relations of two independent solutions? How to define the Wronskian determinant? How to construct a fundamental matrix of the linear QDEs? There is no algorithm to compute the fundamental matrix of the linear QDEs yet. In particular, there is no algorithm to compute the fundamental matrix by using the eigenvalue and eigenvector theory.

In this paper, we try to establish a basic framework for the nonautonomous linear homogenous two-dimensional QDEs taking form of
\[
\left\{
\begin{array}{ccc}
\dot{x}_1& =& a_{11}(t) x_1 +a_{12}(t) x_2
\\
\dot{x}_2 &= & a_{21}(t) x_1 +a_{22}(t) x_2 ,
\end{array}
\right.
\]
where all the coefficients $a_{11},a_{12},a_{21},a_{22}$ are quaternionic functions.
A systematic basic theory on the solution to the higher dimensional linear QDEs is given in this paper. We estabilsh the superposition principle and the structure of the general solutions. Moreover, we provide two algorithms to compute the fundamental matrix of the linear QDEs. For the sake of  understanding, we focus on the 2-dimensional system. However, in fact, our basic results can be easily extended to any arbitrary $n$-dimensional QDEs.

In the following, we will show that there are four profound differences between QDEs and ODEs.

1. On the non-commutativity of the quaternion algebra, the algebraic structure of the solutions to QDEs is {\bf not a linear vector space}. It is actually a {\bf right free module}. Consequently, to determine if two solutions of the QDEs are linearly independent or dependent, we have to define the concept of left dependence (independence) and right dependence (independence). 
 While for classical ODEs, left dependence and right dependence are the same.

2. {\em Wronskian} of ODEs is defined by Caley determinant. However, since Caley determinant depends on the expansion of $i$-th row and $j$-th column of quaternion, different expansions can lead to different results. {\em Wronskian} of ODEs can not be extended to QDEs. It is necessary to define the {\em Wronskian} of QDEs by a novel method (see Section 4). We use double determinants $\rm {ddet}M=\displaystyle\frac{1}{2}\rm{rddet}(MM^{+})$ or $\displaystyle\frac{1}{2}\rm{rddet}(M^{+}M)$ to define {\em Wronskian} of QDEs.

3. Liouville formula for QDEs and ODEs are different.

4. For QDEs, it is necessary to treat the eigenvalue problems with left and right, separately. This is a large-difference between QDEs and ODEs.

 \subsection{Outline of the paper}

 In next section, quaternion algebra is introduced. In section 3, existence and uniqueness of the solutions to QDEs are proved. Section 4 devotes to defining the Wronskian determinant for QDEs and proving the Liouville formulas. Superposition principle and structure of the general solution are given. In section 5, fundamental matrix and solution to linear QDEs are given. Section 6 presents two algorithms for computing fundamental matrix.
Some examples are given to show the feasibility of the obtained
algorithms. Finally, a conclusion and discussion ends the paper.

\section{\bf Quaternion algebra}

The quaternions were discovered in 1843 by Sir William Rowan Hamilton \cite{Hamilton}. We adopt the standard notation in \cite{Chen1,Hamilton}. We denote as usual a quaternion $q=(q_0, q_1, q_2, q_3)^T\in \mathbb{R}^4$ by
\[
q= q_0 + q_1 \mathbf{i} + q_2 \mathbf{j} + q_3 \mathbf{k},
\]
where $q_0, q_1, q_2, q_3$ are real numbers and $\mathbf{i},\mathbf{j},\mathbf{k}$ symbols satisfying the multiplication table formed by
\[
\mathbf{i}^2=\mathbf{j}^2=\mathbf{k}^2=\mathbf{ijk}=-1.
\] We denote the set of quaternion by $ \mathbb{H}$.

We can define for the quaternion $q$ and $h=h_0 + h_1 \mathbf{i} + h_2 \mathbf{j} + h_3 \mathbf{k}$, the {\em inner product}, and the {\em modulus}, respectively, by
\[
< q, h>= q_0 h_0+ q_1 h_1 + q_2 h_2 +q_3 h_3,
\]
\[
|q|= \sqrt {< q, h>}=\sqrt {q_0^2 + q_1^2 + q_2^2 + q_3^2 }.
\]
\[
\Re q= q_0,\,\,\,\, \Im q= q_1 \mathbf{i} + q_2 \mathbf{j} + q_3 \mathbf{k}.
\]
And the conjugate is
\[
\overline {q}= q_0 - q_1 \mathbf{i} - q_2 \mathbf{j} - q_3 \mathbf{k}.
\]
Thus, we have
\[
\overline {q h}= \overline{h} \overline {q},
\,\,\,\,\, and \,\,\,\,
q^{-1}=\frac{\overline{q}} {|q|^2},\,\,\,\,q\neq0.
\]
For any $q,p\in \mathbb{H}$, we have
\[
|q+p|\leq|q|+|p|,\,\,\,|q\cdot p|\leq|q|\cdot|p|.
\]
Let $\psi: \mathbb{R}  \rightarrow \mathbb{H}$ be a quaternion-valued function defined on $\mathbb{R}$ ($t\in \mathbb{R}$ is a real variable). We denote the set of such quaternion-valued functions by $\mathbb{H}\otimes \mathbb{R}$.
Denote
$
\mathbb{H}^2= \{\Psi \big| \Psi =  \left ( \begin{array}{l} \psi_1
\\
\psi_2
\end{array}
\right )
\}$.
 Then two dimensional quaternion-valued functions defined on real domain, $\Psi (t)=  \left ( \begin{array}{l} \psi_1(t)
\\
\psi_2(t)
\end{array}
\right ) \in   \mathbb{H}^2 \otimes \mathbb{R},$
where $\psi_i (t)=   \psi_{i0} (t) + \psi_{i1}  (t) \mathbf{i}
    + \psi_{i2}(t) \mathbf{j} +   \psi_{i3} (t) \mathbf{k} $.
And the first derivative of two dimensional quaternionic functions with respect to the real variable $t$ will be denoted by
\[
\dot \Psi (t) =  \left ( \begin{array}{l} \dot {\psi}_1 (t)
\\
\dot {\psi}_2 (t)
\end{array}
\right ),\,\,\,
    \dot {\psi}_i:=\frac { {\rm d} \psi_i } { {\rm d}  t } =\frac { {\rm d} \psi_{i0} } { {\rm d}  t } +\frac { {\rm d} \psi_{i1} } { {\rm d}  t } \mathbf{i}
    +\frac { {\rm d} \psi_{i2} } { {\rm d}  t } \mathbf{j} +\frac { {\rm d} \psi_{i2} } { {\rm d}  t } \mathbf{k} .
\]
 For any $A(t)=(a_{ij}(t))_{n\times n}\in \mathbb{H}^{n\times n} \otimes \mathbb{R}$,
 \[
 \frac{ {\rm d} A(t)} {{\rm d} t}=\left(\frac{ {\rm d} a_{ij}(t)}{{\rm d} t} \right).
 \]
 For any quaternion-valued functions $f(t)$ defined on $\mathbb{R}$, $f(t)=f_{0}(t)+f_{1}(t)\mathbf{i}+f_{2}(t)\mathbf{j}+f_{3}(t)\mathbf{k}$ on the interval $[a,b]\subset \mathbb{R}$, if $f_{l}(t)$ $(l=0,1,2,3)$ is integrable, then $f(t)$ is integrable and
\[
\int^{b}_{a}f(t)dt=\int^{b}_{a}f_{0}(t)dt+\int^{b}_{a}f_{1}(t)dt\mathbf{i}+\int^{b}_{a}f_{2}(t)dt\mathbf{j}+\int^{b}_{a}f_{3}(t)dt\mathbf{k}, (a\leq t \leq b).
\]
Moreover,
 \[
  \int_{a}^{b} A(t)dt=\left(  \int_{a}^{b} a_{ij}(t)dt \right).
 \]
We define the norms for the quaternionic vector and matrix.
 For any $A=(a_{ij})_{n\times n}\in \mathbb{H}^{n\times n}$ and ${x}=(x_{1},x_{2},\cdots,x_{n})^{T}\in \mathbb{H}^{n}$, we define
 \[
 \|A\|=\sum\limits^{n}_{i,j=1}|a_{ij}|,\,\,\,\| {x}\|=\sum\limits^{n}_{i=1}|x_{i}|.
 \]
For $A,B\in \mathbb{H}^{2\times 2}$ and $ {x}, {y}\in \mathbb{H}^{2}$, we see that
\[
\begin{array}{ccc}
\|AB\|&\leq&\|A\|\cdot\|B\|,
\\
\|A {x}\|&\leq&\|A\| \cdot\| {x}\|,
\\
\|A+B\|&\leq&\|A\|+\|B\|,
\\
\| {x}+ {y}\|&\leq&\| {x}\|+ \| {y}\|.
\end{array}
\]

\section{Existence and Uniqueness of Solution to QDEs}

Consider the linear system
\begin{equation}
\dot \Psi(t) = A(t)\Psi(t)\,\,\,\, or \,\,\,\,\,\,
\left ( \begin{array}{l} \dot {\psi}_1 (t)
\\
\dot {\psi}_2 (t)
\end{array}
\right )= \left ( \begin{array}{ll}   {a} _{11}(t)  & {a} _{12}(t)
 \\
 {a} _{21}(t)  & {a} _{22}(t)
\end{array}
\right ) \left ( \begin{array}{l}   {\psi}_1 (t)
\\
  {\psi}_2 (t)
\end{array}
\right ),
\label{2.1}
\end{equation}
 where $  \Psi  (t) \in   \mathbb{H}^2 \otimes \mathbb{R}$, $A(t)  \in   \mathbb{H}^{2\times 2} \otimes \mathbb{R}$ is $2 \times 2$ continuous 
 quaternion-valued function matrix on the interval $I=[a,b],\, I\subset \mathbb{R}$.

System (\ref{2.1}) is associated with the following initial value problem
\begin{equation}
\Psi (t_0) =  \xi, \,\,\,\,\,\,  \xi  \in \mathbb{H}^2 .
\label{2.3}
 \end{equation}
\indent To discuss the existence and uniqueness of the solutions to QDEs (\ref{2.1}), we need introduce the limit of the quaternion-valued sequence and convergence of the quaternion-valued sequence. Let $\{x^{k}\}\in \mathbb{H}^n$ be a sequence of quaternionic vector, $x^k= (x^k_1,x^k_2,\cdots,x^k_n)^T$, where $x^k_{i}=x^k_{i0}+ x^k_{i1}\mathbf{i} + x^k_{i2}\mathbf{j} + x^k_{i3}k$, $x^k_{il}, (i = 1,2,\cdots,n;l = 0,1,2,3)$ are real numbers. If the sequences of
real numbers $\{x^k_{il}\}$ are convergent, then the sequence of quaternonic vectors $\{x^{k}\}$ is said to be
convergent. And we define
 \[
  \lim\limits_{k\rightarrow \infty}x^k_{i}=\lim\limits_{k\rightarrow \infty}x_{i0}^{k}+\lim\limits_{k\rightarrow \infty}x_{i1}^{k}\mathbf{i}+\lim\limits_{k\rightarrow \infty}x_{i2}^{k}\mathbf{j}+\lim\limits_{k\rightarrow \infty}x_{i3}^{k}\mathbf{k}.
  \]
Let $\{x^{k}(t)\}\in \mathbb{H}^n \otimes \mathbb{R}$ be a sequence of quaternion-valued vector functions defined on $\mathbb{R}$.
$x^k(t)= (x^k_1(t),x^k_2(t),\cdots,x^k_n(t))^T$, where $x^k_{i}(t)=x^k_{i0}(t)+ x^k_{i1}(t)\mathbf{i} + x^k_{i2}(t)\mathbf{j} + x^k_{i3}(t)k$, $x^k_{il}(t)\in \mathbb{R}^n \otimes \mathbb{R}, (i = 1,2,\cdots,n;l = 0,1,2,3)$ are real functions. If the sequences of real-valued functions $\{x^k_{il}(t)\}$
are convergent (uniform convergent), then the sequence of quaternion-valued vector functions
$\{x^{k}(t)\}$ is said to be convergent (uniform convergent). Moreover,
 \[
  \lim\limits_{k\rightarrow \infty}x^k_{i}(t)=\lim\limits_{k\rightarrow \infty}x_{i0}^{k}(t)+\lim\limits_{k\rightarrow \infty}x_{i1}^{k}(t)\mathbf{i}+\lim\limits_{k\rightarrow \infty}x_{i2}^{k}(t)\mathbf{j}+\lim\limits_{k\rightarrow \infty}x_{i3}^{k}(t)\mathbf{k}.
  \]
If the sequence of vector functions $S_n(t) =
\sum\limits ^{n}_{k=1} {x}^{k}(t)$ is convergent (uniform convergent), then
the series
$\sum\limits ^{\infty}_{k=1} {x}^{k}(t)$ is said to be convergent (uniform convergent). It is not difficult to prove
that if there exist positive real constants $\delta_k$ such that
 \[
  \| {x}^{k}(t)\|\leq \delta_{k}, \,\,(a\leq t \leq b)
  \]
and the real series
$\sum\limits ^{\infty}_{k=1}\delta_{k}$ is convergent, then the quaternion-valued series
 $\sum\limits ^{\infty}_{k=1} {x}^{k}(t)$  is absolutely
convergent on the interval $[a,b]$. Certainly, it is uniform convergent on the interval $[a,b]$.
Moreover, if the sequence of quaternion-valued vector functions $\{x^k(t)\}$ is uniform convergent
on the interval $[a,b]$, then
\[
  \lim_{k\rightarrow \infty}\int^{b}_{a} {x}^{k}(t)dt=\int^{b}_{a}\lim_{k\rightarrow \infty} {x}^{k}(t)dt.
  \]
Based on the above concepts on the convergence of quaternion-valued series, by similar arguments to Theorem 3.1 in Zhang et al \cite{ZWN}, we have
\begin{theorem}\label{Th3.1}
  The initial value problem (\ref{2.1})-(\ref{2.3}) has exactly a unique solution.
\end{theorem}

\section{ Wronskian and structure of general solution to QDEs }

To study the superposition principle and structure of the general solution to QDEs, we should introduce the concept of right- module.
On the non-commutativity of the quaternion algebra, we will see that the algebraic structure of the solutions to QDEs is {\bf not a linear vector space}. It is actually a right- module. This is the largest difference between QDEs and ODEs.

Now we introduce some definitions and known results on groups, rings, modules in the abstract algebraic theory (e.g. \cite{abstract}).
An {\em abelian group} is a set, $\mathbb{A}$, together with an operation $\cdot$ that combines any two
elements $a$ and $b$ to form another element denoted $a \cdot b$. The symbol $\cdot$ is a general
placeholder for a concretely given operation. The set and operation, $(\mathbb{A}, \cdot)$, must satisfy five requirements known as the abelian group
axioms:
\begin{description}

\item[(i)] {\em  Closure:} For all $a, b$ in $\mathbb{A}$, the result of the operation $a \cdot b$ is also in $\mathbb{A}$.

\item[(ii)] {\em Associativity:} For all $a, b$ and $c$ in $\mathbb{A}$, the equation $(a \cdot b) \cdot c = a\cdot (b \cdot c)$ holds.

\item[(iii)] {  \em Identity element:} There exists an element $e$ in $\mathbb{A}$, such that for all elements a in $\mathbb{A}$, the equation
$e \cdot a = a \cdot e = a$ holds.

\item[(iv)] { \em Inverse element:} For each a in $\mathbb{A}$, there exists an element b in $\mathbb{A}$ such that $a \cdot b = b \cdot a = e$,
where $e$ is the identity element.

\item[(v)] { \em Commutativity:} For all $a, b$ in $\mathbb{A}$, $a \cdot b = b \cdot a$.
\end{description}

More compactly, an abelian group is a commutative group.

\smallskip

A {\em ring} is a set $(\cal{R},+, \cdot)$ equipped with binary operations $+$ and $\cdot$ satisfying the
following three sets of axioms, called the ring axioms:
\begin{description}

\item[(i)] $\cal{R}$ is an abelian group under addition, meaning that

$\bullet $   $(a + b) + c = a + (b + c)$ for all $a, b, c$ in $\cal{R}$ ($+$ is associative).

$\bullet $   $a + b = b + a$ for all $a, b$ in $\cal{R}$ ($+$ is commutative).

$\bullet $   There is an element $0$ in $\cal{R}$ such that $a + 0 = a$ for all a in $\cal{R}$ ($0$ is the
additive identity).

$\bullet $  For each $a$ in $\cal{R}$ there exists .a in R such that $a + (-a) = 0$ ($-a$ is the
additive inverse of $a$).

\item[(ii)] $\cal{R}$ is a monoid under multiplication, meaning that:

$\bullet $  $(a \cdot b) \cdot c = a \cdot (b \cdot c)$ for all $a, b, c$ in $\cal{R}$ ($\cdot$ is associative).

$\bullet $ There is an element $1$ in $\cal{R}$ such that $a \cdot 1 = a$ and $1 \cdot a = a$ for all $a$ in
$\cal{R}$ ($1$ is the multiplicative identity).

\item[(iii)] Multiplication is distributive with respect to addition:

$\bullet $  $a \cdot (b + c) = (a \cdot b) + (a \cdot c)$ for all $a, b, c$ in $\cal{R}$ (left
distributivity).

$\bullet $  $(b + c) \cdot  a = (b \cdot a) + (c \cdot a)$ for all $a, b, c$ in $\cal{R}$ (right
distributivity).
\end{description}

 $(\cal{F},+, \cdot)$ is said to be a field if $(\cal{F},+, \cdot)$ is a ring and ${\cal{F}}^{*}={\cal{F}}/\{0\}$ with respect to multiplication $\cdot$ is a commutative group.  The most commonly used fields are the field of real numbers, the field of complex numbers.
A ring in which division is possible but commutativity is not assumed (such as the quaternions) is called a division ring. In fact, a field is a nonzero commutative division ring.

\begin{remark}\label{Rem4.1}
Different from real or complex numbers, quaternion is a division ring, not a field due to the non-commutativity. Thus, the vector space $\mathbb{H}^2$ over quaternion is not a linear vector space. We need introduce the concept of {\bf right- module}.
\end{remark}

Module over a ring $\cal{R}$ is a generalization of abelian group. You may
view an $\cal{R}$-mod as a ``vector space over $\cal{R}$". A {\em right $\cal{R}$-module} is an additive abelian group $\mathbb{A}$
together with a function $\mathbb{A}\times \cal{R} \rightarrow \mathbb{A}$ (by $(a, r) \rightarrow a\cdot r$ ) such that for all
$a, b \in  \mathbb{A}$ and $r,s \in \cal{R}$:
\begin{description}
\item[(i)] $(a + b)\cdot r =  a \cdot r +  b \cdot r$.
\item[(ii)] $a \cdot (r + s) = a \cdot r +  a \cdot  s$.
\item[(iii)]  $(a\cdot r)\cdot s = a\cdot (rs)$.
\end{description}
Similarly, we can define the left $\cal{R}$-module.

The multiplication $xr$, where $x\in\mathbb{A}$ and $r\in\cal{R}$, is called the right action of $\cal{R}$ on the abelian group $\mathbb{A}$. We will write $\mathbb{A}_{\cal{R}}^R$ to indicate that $\mathbb{A}$ is a right $\cal{R}$-module (over ring $\cal{R}$). Similarly, we can define the left $\cal{R}$-module. And we denote the left $\cal{R}$-module by  $\mathbb{A}_{\cal{R}}^L$ (over ring $\cal{R}$).

\textbf{Right Submodules}. Let $\mathbb{A}_{\cal{R}}^R$ be a right $\cal{R}$-module, a subset $\mathbb{N}^{R}\subset\mathbb{A}_{\cal{R}}^R$ is called a right $\cal{R}$-submodule, if the following conditions are satisfied:
\begin{description}
\item[(i)] $\mathbb{N}^{R}$ is a subgroup of the (additive, abelian) group $\mathbb{A}^{R}$.
\item[(ii)] $xr$ is in $\mathbb{N}^{R}$ for all $r\in \cal{R}$ and $x\in\mathbb{N}^{R}$.
\end{description}
Similarly, we can define the concept left submodules.

\textbf{Direct sums}. If $\mathbb{N}_{1}, \cdots, \mathbb{N}_{m}$ are submodules of a module $\mathbb{A}_{\cal{R}}$ over a ring $\cal{R}$, their sum $\mathbb{N}_{1}+\cdots+\mathbb{N}_{m}$ is defined to be the set of all sums of elements of the $\mathbb{N}_{i}$, namely,
\[
\mathbb{N}_{1}+\cdots+\mathbb{N}_{m}=\{x_{1}+\cdots+x_{m}|x_{i}\in\mathbb{N}_{i}\,\, for \,\, each \,\, i\}
\]
this is a submodule of $\mathbb{A}$. We say $\mathbb{N}_{1}+\cdots+\mathbb{N}_{m}$ is a direct sum if this representation is unique for each $x\in\mathbb{N}_{1}+\cdots+\mathbb{N}_{m}$; that is, given two representations
\[
x_{1}+\cdots+x_{m}=x=y_{1}+\cdots+y_{m}
\]
where $x_{i}\in\mathbb{N}_{i}$ and $y_{i}\in\mathbb{N}_{i}$ for each $i$, then necessarily $x_{i}=y_{i}$ for each $i$.

\begin{lemma}\label{Lemma4.0}
 The following conditions are equivalent for a set $\mathbb{N}_{1}, \cdots, \mathbb{N}_{m}$ of submodules of a module $\mathbb{A}$:
\begin{description}
\item[(i)] The sum $\mathbb{N}_{1}+\cdots+\mathbb{N}_{m}$ is direct.
\item[(ii)] $(\sum\limits_{i\neq k}\mathbb{N}_{i})\cap \mathbb{N}_{k}=0$ for each $k$.
\item[(iii)] $\mathbb{N}_{1}+\cdots+\mathbb{N}_{k-1}\cap \mathbb{N}_{k}=0$ for each $k\geq2$.
\item[(iv)] If $x_{1}+\cdots+x_{m}=0$ where $x_{i}\in\mathbb{N}_{i}$ for each $i$, then each $x_{i}=0$
\end{description}
If $\mathbb{N}_{1}+\cdots+\mathbb{N}_{m}$ is direct sum, we denote it by $\mathbb{N}_{1}\oplus\cdots\oplus\mathbb{N}_{m}$.
\end{lemma}

Given elements $x_{1}, \cdots, x_{k}$ in a right- module $\mathbb{A}_{\cal{R}}^{R}$ over ring $\cal{R}$, a sum $x_{1}r_{1}+ \cdots+x_{k}r_{k}$, $r_{i}\in{\cal{R}}$, is called right- linear combination of the $x_{i}$, with coefficients $r_{i}$. The sum of all the right- submodules $x_{i}{\cal{R}},\,(i=1,2,\cdots,k)$ is a right- submodule
\[
x_1\mathcal R+\cdots+ x_k\mathcal R=\{x_1 r_1+\cdots+x_k r_k |r_i \in \mathcal R\}
\]
consisting of all such right linear combinations. If $\mathbb{A}_{\cal{R}}^{R}=x_1\mathcal R+\cdots+x_k\mathcal R$, we say that $\{x_1, \cdots, x_k\}$ is a generating set for in a right- module $\mathbb{A}_{\cal{R}}^{R}$. Similarly, the sum of all the left submodules ${\cal{R}}x_{i},\,(i=1,2,\cdots,k)$ is a left submodule
\[
{\mathcal R} x_1 + \cdots + {\mathcal R}  x_k =\{ r_1 x_1+\cdots+r_k x_k |r_i \in \mathcal R\}
\]
consisting of all such right linear combinations. If $\mathbb{A}_{\cal{R}}^{L}={\mathcal R} x_1 + \cdots+{\mathcal R}  x_k $, we say that $\{x_1, \cdots, x_k\}$ is a generating set for in a left module $\mathbb{A}_{\cal{R}}^{L}$.

The elements $x_{1}, \cdots, x_{k}\in \mathbb{A}_{\cal{R}}^{R}$ are called independent if
\[
x_1 r_1+ \cdots+x_k r_k=0,r_i\in\mathcal R~\text{implies that}~r_1=\cdots=r_k=0.
\]
For the sake of convenience, we call it {\em right} independence.
 A subset $\{x_{1}, \cdots, x_{k}\}$ is called a basis of a right module $\mathbb{A}_{\cal{R}}^{R}$ if it is right independent and generates  $\mathbb{A}_{\cal{R}}^{R}$. Similarly, the elements $x_{1}, \cdots, x_{k}\in \mathbb{A}_{\cal{R}}^{L}$ are called {\em left independent} if
\[
r_1 x_1+ \cdots+r_k x_k=0,r_i\in\mathcal R~\text{implies that}~r_1=\cdots=r_k=0.
\]
 A subset $\{x_{1}, \cdots, x_{k}\}$ is called a basis of a left module $\mathbb{A}_{\cal{R}}^{L}$ if it is left independent and generates  $\mathbb{A}_{\cal{R}}^{L}$. A module that has a finite basis is called a {\bf free module}.

From the above discussion, we know that the quaternionic vector space $\mathbb{H}^{n}$ over the division ring $\mathbb{H}$ is a right- or left- module.
 Thus, the quaternionic vectors $x_1,x_2,\cdots,x_n$, $x_i\in \mathbb{H}^{n}$ are right independent if
\[
x_1 r_1+ \cdots+x_n r_n=0,r_i\in\mathbb{H}~\text{implies that}~r_1=\cdots=r_k=0.
\]
On the contrary, if there is a nonzero $r_i,\i=1,2,\cdots,n$ such that $x_1 r_1+ \cdots+x_n r_n=0$, then $x_1,x_2,\cdots,x_n$ are rihgt dependent. Similarly, we can define the left independence or dependence for the quaternionic vector $\mathbb{H}^{n}$ over the division ring $\mathbb{H}$.

For the sake of easily understanding, we focus on the two dimensional system (\ref{2.1}).

  \begin{definition}\label{Def4.1}
  For two quaternion-valued functions $x_1(t)$ and $x_2(t)$ defined on the real interval $I$, if there are two quaternionic constants $q_1,q_2\in \mathbb{H}$  (not both zero) such that
\[
x_1(t) q_1 + x_2(t) q_2 =0,\,\,  \,\,\,\, for \,\,\,\,any \,\,\,\, t\in I,
\]
then $x_1(t)$ and $x_2(t)$ are said to be {\em right} linearly dependent on $I$. On the contrary, $x_1(t)$ and $x_2(t)$ are said to be right linearly {\em independent} on $I$ if the algebric equation
\[
x_1(t) q_1 + x_2(t) q_2 =0,\,\,  \,\,\,\, for \,\,\,\,any \,\,\,\, t\in I,
\]
can only be satisfied by
\[
q_1=q_2=0.
\]
If there are two quaternionic constants $q_1,q_2\in \mathbb{H}$  (not both zero) such that
\[
  q_1 x_1(t)  + q_2 x_2(t)  =0,\,\,\,\, for \,\,\,\,any \,\,\,\, t\in I,
\]
then  $x_1(t)$ and $x_2(t)$ are said to be {\em left} linearly dependent on $I$. $x_1(t)$ and $x_2(t)$ are said to be left linearly {\em independent} on $I$ if the algebric equation
\[
  q_1 x_1(t)  + q_2 x_2(t)  =0,\,\,\,\, for \,\,\,\,any \,\,\,\, t\in I,
\]
can only be satisfied by
$
q_1=q_2=0.
$
\end{definition}

For sake of convenience, in this paper, we adopt Cayley determinant over the quaternion \cite{Cayley}.  For a $2\times 2$ determinant we could use ${a}_{11} {a} _{22} - {a} _{12} {a} _{21}$ (expanding along the first row). That is,
\[
{\rm {rdet}}
 \left ( \begin{array}{ll}
  {a} _{11}  & {a} _{12}
 \\
 {a} _{21}  & {a} _{22}
\end{array}
\right )
=
 {a}_{11} {a} _{22} - {a} _{12} {a} _{21}.
 \]
Certainly, one can also use ${a}_{11} {a} _{22} - {a} _{21} {a} _{12}$ (expanding along the first column) as follows.
\[
{\rm {cdet}}
 \left ( \begin{array}{ll}
  {a} _{11}  & {a} _{12}
 \\
 {a} _{21}  & {a} _{22}
\end{array}
\right )
=
 {a}_{11} {a} _{22} - {a} _{21} {a} _{12}.
 \]
 There are other definitions (expanding along the second column) or (expanding along the second row).
 Without loss of generality, we employ ${\rm rdet}$ to proceed our study in this paper.
In fact, one can prove the main results in this paper by replacing the Cayley determinant by other definitions.

To study the solution of Eq. (\ref{2.1}), firstly, we should define the concept of {\em Wronskian} for QDEs.
Consider any two solutions of (\ref{2.1}),
\[
x_1(t)
=
\left ( \begin{array}{l}   {x} _{11}(t)
\\
  {x}_{21} (t)
\end{array}
\right ),\,\,\,\,
x_2(t)
=
\left ( \begin{array}{l}   {x} _{12}(t)
\\
  {x}_{22} (t)
\end{array}
\right )
\in
\mathbb{H}^2 \otimes \mathbb{R}.
\]
If we define {\em Wronskian} as standard consideration of ODEs,
\begin{equation}
 W_{ODE}(t) =
{{\rm {rdet}}}
 \left ( \begin{array}{ll}
  {x} _{11}(t)  & {x} _{12}(t)
 \\
 {x} _{21}(t)  & {x} _{22}(t)
\end{array}
\right )
=
 {x} _{11}(t) {x} _{22}(t) - {x} _{12}(t)
 {x} _{21}(t).
 \label{W-ODE}
 \end{equation}
This Wronskian definition cannot be extended to QDEs. In fact, let us consider two right
linearly dependent solutions of (\ref{2.1}),
\begin{equation}
x_1(t)
= x_2(t) \eta,
\,\,\,\, or \,\,\,\,
\left ( \begin{array}{l}   {x} _{11}(t)
\\
  {x}_{21} (t)
\end{array}
\right )
=
\left ( \begin{array}{l}   {x} _{12}(t)
\\
  {x}_{22} (t)
\end{array}
\right )\eta,\,\,\,\,
x_1(t)
,x_2(t)
\in
\mathbb{H}^2 \otimes \mathbb{R},\,\,\, \eta \in \mathbb{H},
\label{LI}
 \end{equation}
which implies
\[
x_{12}(t)= x_{11}(t) \eta^{-1},
\,\,\, and \,\,\,
x_{21}(t)= x_{22}(t) \eta.
\]
Substituting these inequalities into (\ref{W-ODE}), we have
\[
 W_{ODE}(t) = {x} _{11}(t) {x} _{22}(t) - {x} _{12}(t)
 {x} _{21}(t)= {x} _{11}(t) {x} _{22}(t) - {x} _{11}(t)\eta^{-1}
 {x} _{22}(t)\eta \neq 0.
\]
The reason lies in $\eta$ is a quaternion, which can not communicate. Therefore, we should change the standard definition of {\em Wronskian}. Now we give the concept of {\em Wronskian} for QDEs as follows.

\begin{definition}\label{Def4.2}
 Let $x_1(t),x_2(t)$ be two solutions of Eq. (\ref{2.1}). Denote
\[
M(t)
=
\left ( \begin{array}{ll}
  {x} _{11}(t)  & {x} _{12}(t)
 \\
 {x} _{21}(t)  & {x} _{22}(t)
\end{array}
\right ).
\]
The {\em Wronskian} of QDEs is defined by
\begin{equation}
 W_{QDE}(t) =\rm{ddet} M(t): = \frac{1}{2}{\rm {rdet}} \big(M(t) M^+(t) \big)
 \label{W}
 \end{equation}
where $\rm{ddet} M:=\displaystyle {\rm {rdet}} \big(M  M^+  \big)$ is called double determinant, $M^{+}$ is the conjugate transpose of $M$, that is,
\[
M^{+}(t)
=
\left ( \begin{array}{ll}
 \overline {x} _{11}(t)  &  \overline {x} _{21}(t)
 \\
\overline {x} _{12}(t)  & \overline {x} _{22}(t)
\end{array}
\right ).
\]
\end{definition}

\begin{remark}\label{Rem4.2}
 From Definition \ref{Def4.2}, we have
\begin{equation}
\begin{array}{lll}
 W_{QDE}(t)
  &=&
  \displaystyle \{ |x_{11}(t)|^2|x_{22}(t)|^2  + |x_{12}(t)|^2|x_{21}(t)|^2
 \\&&
 \hspace{0.5cm}
 - x_{12}(t)   \overline{x}_{22}(t) x_{21}(t) \overline{x}_{11}(t)
 -  x_{11}(t) \overline{x}_{21} (t)x_{22}(t)  \overline{x}_{12} (t) \}.
 \end{array}
\label{W2}
 \end{equation}
We see that
\[
\overline{  x_{12}(t)   \overline{x}_{22}(t) x_{21}(t) \overline{x}_{11}(t) } = x_{11}(t) \overline{x}_{21} (t)x_{22}(t)  \overline{x}_{12} (t),
\]
which, combining with (\ref{W2}), implies that
$W_{QDE}(t)$ is a real number.
\end{remark}

  For quaternionic matrix, it should be pointed out that
\[
{\rm {rdet}} \big(M(t) M^+(t) \big)\neq {\rm {rdet}} M(t) \cdot {\rm {rdet}}   M^+(t)  .
\]

\begin{proposition}\label{Pro4.1}
  If $x_1(t)$ and $x_2(t)$ are right 
   linearly dependent on $I$ then $W_{QDE}(t)=0$.
  \end{proposition}

{\bf PROOF.}   If $x_1(t)$ and $x_2(t)$ are right linearly dependent on $I$ then Eq. (\ref{LI}) holds. Then by (\ref{W2}),
\[
\begin{array}{lll}
 W_{QDE}(t)
  &=&
  \displaystyle

  \{
  |x_{12}(t)|^2|\eta|^2 |x_{22}(t)|^2+ |x_{11}(t)|^2|\eta|^2|x_{22}(t)|^2
  \\
  &&
  \displaystyle
-  x_{12}(t)   \overline{x}_{22} (t)   x_{22}(t) |\eta|^2  \overline{x}_{12}(t)
-  x_{12}(t)  |\eta|^2  \overline{x}_{22}(t)    x_{22}(t)  \overline{x}_{12} (t)
\} =0.
\end{array}
\]
That is, $ W_{QDE}(t)=0$.



\begin{theorem}\label{Th4.1}
 (Liouville formula) The Wronskian  $W_{QDE}(t)$ of Eq. (\ref{2.1}) satisfies the following quaternionic Liouville formula.
 \begin{equation}
W_{QDE}(t)= \exp\left( \int_{t_0}^{t}{\rm tr} [A(s)+A^{+}(s)] ds\right) W_{QDE}(t_0),
\label{Liou}
 \end{equation} where ${\rm tr} A(t)$ is the trace of the coefficient matrix $A(t)$, i.e. ${\rm tr} A(t)=a_{11}(t)+a_{22}(t)$.
   Moreover, if $W_{QDE}(t)=0$ at some $t_0$ in $I$ then $W_{QDE}(t)=0$ on $I$.
 \end{theorem}


{\bf PROOF.}   Let us consider Eq. (\ref{W}). By calculating the first derivative of the left-hand-side and right-hand-side terms, by using (\ref{W2}), we obtain
\[
\begin{array}{lll}
&&
\displaystyle\frac{d}{dt}  W_{QDE}(t)
\\
&=&
 \displaystyle  \frac{d}{dt} {\rm {rdet}}  \left( M(t) M^{+}(t) \right)
  \\
  &=&
    \displaystyle    \frac{d}{dt}\Big[ x_{11}(t) \overline{x}_{11}(t)  x_{22}(t) \overline{x}_{22}(t)  +  x_{12}(t) \overline{x}_{12}(t)  x_{21}(t)  \overline{x}_{21}(t)
    \\
    &&
    \displaystyle
      \hspace{0.6cm}
 - x_{12}(t)   \overline{x}_{22}(t) x_{21}(t) \overline{x}_{11}(t)
 -  x_{11}(t) \overline{x}_{21} (t)x_{22}(t)  \overline{x}_{12} (t) \Big]
   \\
  &=&
    \displaystyle  \frac{d}{dt} x_{11}(t) \overline{x}_{11}(t)  x_{22}(t) \overline{x}_{22}(t)
     +  x_{11}(t) \displaystyle\frac{d}{dt} \overline{x}_{11}(t)  x_{22}(t) \overline{x}_{22}(t)
    \\
    &&
    + x_{11}(t) \overline{x}_{11}(t) \displaystyle\frac{d}{dt} x_{22}(t) \overline{x}_{22}(t)
    + x_{11}(t) \overline{x}_{11}(t)  x_{22}(t) \displaystyle\frac{d}{dt} \overline{x}_{22}(t)
      \\
&&
     +\displaystyle\frac{d}{dt}
     x_{12}(t) \overline{x}_{12}(t)  x_{21}(t)  \overline{x}_{21}(t)
      +
     x_{12}(t) \displaystyle\frac{d}{dt} \overline{x}_{12}(t)  x_{21}(t)  \overline{x}_{21}(t)
     \\
     &&
     +
     x_{12}(t) \overline{x}_{12}(t) \displaystyle\frac{d}{dt} x_{21}(t)  \overline{x}_{21}(t)
     +
     x_{12}(t) \overline{x}_{12}(t)  x_{21}(t) \displaystyle\frac{d}{dt} \overline{x}_{21}(t)
    \\
    &&
  -  \displaystyle\frac{d}{dt}
x_{12}(t)   \overline{x}_{22}(t) x_{21}(t) \overline{x}_{11}(t)
 - x_{12}(t) \displaystyle\frac{d}{dt}  \overline{x}_{22}(t) x_{21}(t) \overline{x}_{11}(t)
 \\
 &&
 - x_{12}(t)   \overline{x}_{22}(t) \displaystyle\frac{d}{dt}x_{21}(t) \overline{x}_{11}(t)
 - x_{12}(t)   \overline{x}_{22}(t) x_{21}(t) \displaystyle\frac{d}{dt}\overline{x}_{11}(t)
 \\
 &&
 - \displaystyle\frac{d}{dt}
x_{11}(t) \overline{x}_{21} (t)x_{22}(t)  \overline{x}_{12} (t)
 -  x_{11}(t) \displaystyle\frac{d}{dt}\overline{x}_{21} (t)x_{22}(t)  \overline{x}_{12} (t)
 \\
 &&
 -  x_{11}(t) \overline{x}_{21} (t) \displaystyle\frac{d}{dt} x_{22}(t)  \overline{x}_{12} (t)
 -  x_{11}(t) \overline{x}_{21} (t)x_{22}(t) \displaystyle\frac{d}{dt} \overline{x}_{12} (t)
\\
  &=&
   \displaystyle  [a_{11}(t) x_{11}(t)+ a_{12}(t)x_{21}(t)] \overline{x}_{11}(t)  x_{22}(t) \overline{x}_{22}(t)
   \\
   &&
     +  x_{11}(t)\displaystyle [\overline{x}_{11}(t) \overline{a}_{11}(t)+ \overline{x}_{21}(t) \overline{a}_{12}(t)]  x_{22}(t) \overline{x}_{22}(t)
    \\
    &&
    + x_{11}(t) \overline{x}_{11}(t)\displaystyle [a_{21}(t) x_{12}(t)+ a_{22}(t)x_{22}(t)]  \overline{x}_{22}(t)
    \\
    &&
    + x_{11}(t) \overline{x}_{11}(t)  x_{22}(t)  \displaystyle [ \overline{x}_{12}(t) \overline{a}_{21}(t)+ \overline{x}_{22}(t) \overline{a}_{22}(t)]
      \\
&&
     +\displaystyle[a_{11}(t) x_{12}(t)+ a_{12}(t)x_{22}(t)]  \overline{x}_{12}(t)  x_{21}(t)  \overline{x}_{21}(t)
       \\
      &&
      +
     x_{12}(t) \displaystyle[\overline{x}_{12}(t) \overline{a}_{11}(t)+ \overline{x}_{22}(t) \overline{a}_{12}(t)]   x_{21}(t)  \overline{x}_{21}(t)
     \\
     &&
     +
     x_{12}(t) \overline{x}_{12}(t) \displaystyle [a_{21}(t) x_{11}(t)+ a_{22}(t)x_{21}(t)]  \overline{x}_{21}(t)
     \\
     &&
     +
     x_{12}(t) \overline{x}_{12}(t)  x_{21}(t) \displaystyle [ \overline{x}_{11}(t)  \overline{a}_{21}(t)+  \overline{x}_{21}(t)  \overline{a}_{22}(t)]
    \\
    &&
  -  \displaystyle[a_{11}(t) x_{12}(t)+ a_{12}(t)x_{22}(t)]    \overline{x}_{22}(t) x_{21}(t) \overline{x}_{11}(t)
  \\
  &&
 - x_{12}(t) \displaystyle [ \overline{x}_{12}(t) \overline{a}_{21}(t)+ \overline{x}_{22}(t) \overline{a}_{22}(t)] x_{21}(t) \overline{x}_{11}(t)
 \\
 &&
 - x_{12}(t)   \overline{x}_{22}(t)  \displaystyle [a_{21}(t) x_{11}(t)+ a_{22}(t)x_{21}(t)] \overline{x}_{11}(t)
 \\
 &&
 - x_{12}(t)   \overline{x}_{22}(t) x_{21}(t) \displaystyle [\overline{x}_{11}(t) \overline{a}_{11}(t)+ \overline{x}_{21}(t) \overline{a}_{12}(t)]
 \\
 &&
 - \displaystyle[a_{11}(t) x_{11}(t)+ a_{12}(t)x_{21}(t)]  \overline{x}_{21} (t)x_{22}(t)  \overline{x}_{12} (t)
 \\
 &&
 -  x_{11}(t)  \displaystyle [ \overline{x}_{11}(t)  \overline{a}_{21}(t)+  \overline{x}_{21}(t)  \overline{a}_{22}(t)] x_{22}(t)  \overline{x}_{12} (t)
 \\
 &&
 -  x_{11}(t) \overline{x}_{21} (t) \displaystyle [a_{21}(t) x_{12}(t)+ a_{22}(t)x_{22}(t)]  \overline{x}_{12} (t)
 \\
 &&
 -  x_{11}(t) \overline{x}_{21} (t)x_{22}(t)  \displaystyle[\overline{x}_{12}(t) \overline{a}_{11}(t)+ \overline{x}_{22}(t) \overline{a}_{12}(t)]
\end{array}
\]
\[
\begin{array}{lll}
  &=&
     \displaystyle \Big\{ {a}_{11}(t) |x_{11}(t)|^2|x_{22}(t)|^2 + a_{12}(t)x_{21}(t) \overline{x}_{11}(t)  |x_{22}(t)|^2
   \\
   &&
     +  |x_{11}(t)|^2 \overline{a}_{11}(t)|x_{22}(t)|^2  +x_{11}(t) \overline{x}_{21}(t) \overline{a}_{12}(t)]  x_{22}(t) \overline{x}_{22}(t)
    \\
    &&
    + |x_{11}(t)|^2  a_{21}(t) x_{12}(t)\overline{x}_{22}(t)  + |x_{11}(t)|^2 a_{22}(t)|x_{22}(t)|^2
    \\
    &&
    + |x_{11}(t)|^2 x_{22}(t)  \overline{x}_{12}(t) \overline{a}_{21}(t)  + |x_{11}(t)|^2  |{x}_{22}(t)|^2 \overline{a}_{22}(t)]
      \\
&&
     + a_{11}(t) |x_{12}(t)|^2 | x_{21}(t)|^2 + a_{12}(t)x_{22}(t)  \overline{x}_{12}(t)  |x_{21}(t)|^2
       \\
      &&
      +
     |x_{12}(t)|^2 \overline{a}_{11}(t) |x_{21}(t)|^2+ x_{12}(t) \overline{x}_{22}(t) \overline{a}_{12}(t)  |x_{21}(t)|^2
     \\
     &&
     +
     |x_{12}(t)|^2  a_{21}(t) x_{11}(t)\overline{x}_{21}(t)  + |x_{12}(t)|^2  a_{22}(t)|x_{21}(t)|^2
     \\
     &&
     +
     |x_{12}(t)|^2  x_{21}(t)  \overline{x}_{11}(t)  \overline{a}_{21}(t) +  |x_{12}(t)|^2 |{x}_{21}(t)|^2  \overline{a}_{22}(t)]
    \\
    &&
  -  a_{11}(t) x_{12}(t) \overline{x}_{22}(t) x_{21}(t) \overline{x}_{11}(t) - a_{12}(t)|x_{22}(t)|^2 x_{21}(t) \overline{x}_{11}(t)
  \\
  &&
 - |x_{12}(t)|^2\overline{a}_{21}(t)x_{21}(t) \overline{x}_{11}(t) - x_{12}(t) \overline{x}_{22}(t) \overline{a}_{22}(t)] x_{21}(t) \overline{x}_{11}(t)
 \\
 &&
 - x_{12}(t)   \overline{x}_{22}(t)  a_{21}(t) |x_{11}(t)|^2  -  x_{12}(t)   \overline{x}_{22}(t)   a_{22}(t)x_{21}(t)] \overline{x}_{11}(t)
 \\
 &&
 - x_{12}(t)   \overline{x}_{22}(t) x_{21}(t) \overline{x}_{11}(t) \overline{a}_{11}(t) - x_{12}(t)   \overline{x}_{22}(t) |{x}_{21}(t)|^2 \overline{a}_{12}(t)]
 \\
 &&
 -  a_{11}(t) x_{11}(t)\overline{x}_{21} (t)x_{22}(t)  \overline{x}_{12} (t) - a_{12}(t)|x_{12}(t)|^2 x_{22}(t)  \overline{x}_{12} (t)
 \\
 &&
 -  |x_{11}(t)|^2   \overline{a}_{21}(t) x_{22}(t)  \overline{x}_{12} (t)  -  x_{11}(t) \overline{x}_{21}(t)  \overline{a}_{22}(t)] x_{22}(t)  \overline{x}_{12} (t)
 \\
 &&
 -  x_{11}(t) \overline{x}_{21} (t) a_{21}(t) |x_{12}(t)|^2 - x_{11}(t) \overline{x}_{21} (t) a_{22}(t)x_{22}(t)]  \overline{x}_{12} (t)
 \\
 &&
 -  x_{11}(t) \overline{x}_{21} (t)x_{22}(t)   \overline{x}_{12}(t) \overline{a}_{11}(t) -   x_{11}(t) \overline{x}_{21} (t) |{x}_{22}(t)|^2 \overline{a}_{12}(t)
 \Big\}
 \end{array}
  \]
  \begin{equation}
\begin{array}{lll}
  &=&
     \displaystyle
     \Big\{
     \big[ {a}_{11}(t) + {a}_{22}(t)  +\overline{a}_{11}(t)  + \overline{a}_{22}(t)  \big] |x_{11}(t)|^2|x_{22}(t)|^2
     \\
     &&
     +
        \big[ {a}_{11}(t) + {a}_{22}(t)  +\overline{a}_{11}(t)  + \overline{a}_{22}(t)  \big] |x_{12}(t)|^2|x_{21}(t)|^2
\\
&&
-{a}_{11}(t)\big[ x_{12}(t) \overline{x}_{22}(t) x_{21}(t) \overline{x}_{11}(t)  + x_{11}(t)\overline{x}_{21} (t)x_{22}(t)  \overline{x}_{12} (t)\big]
\\
&&
-\big[ x_{12}(t) \overline{x}_{22}(t) x_{21}(t) \overline{x}_{11}(t)  + x_{11}(t)\overline{x}_{21} (t)x_{22}(t)  \overline{x}_{12} (t)\big]\overline{a}_{11}(t)
\\
&&
- x_{12}(t) \overline{x}_{22}(t)\big[ {a}_{22}(t)  +\overline{a}_{22}(t) \big]  x_{21}(t) \overline{x}_{11}(t)
\\
&&
- x_{11}(t)\overline{x}_{21} (t)  \big[ {a}_{22}(t)  +\overline{a}_{22}(t) \big] x_{22}(t)  \overline{x}_{12} (t)
 \\
 &&
 +|x_{11}(t)|^2  a_{21}(t) x_{12}(t)\overline{x}_{22}(t) +|x_{11}(t)|^2  {x}_{22}(t) \overline {x}_{12}(t) \overline{a}_{21}(t)
 \\
 &&
 -|x_{11}(t)|^2  \overline{a}_{21}(t) {x}_{22}(t) \overline {x}_{12}(t) - a_{21}(t) x_{12}(t)\overline{x}_{22}(t) |x_{11}(t)|^2
  \\
 &&
 +|x_{12}(t)|^2  a_{21}(t) x_{11}(t)\overline{x}_{21}(t) +|x_{12}(t)|^2  {x}_{21}(t) \overline {x}_{11}(t) \overline{a}_{21}(t)
 \\
 &&
 -|x_{12}(t)|^2  \overline{a}_{21}(t) {x}_{21}(t) \overline {x}_{11}(t) - a_{21}(t) x_{11}(t)\overline{x}_{21}(t) |x_{12}(t)|^2
 \Big\}
  .
       \end{array}
   \label{Deri}
\end{equation}
Set
\[
Q_1=x_{12}(t) \overline{x}_{22}(t) x_{21}(t) \overline{x}_{11}(t)  + x_{11}(t)\overline{x}_{21} (t)x_{22}(t)  \overline{x}_{12} (t)
\]
It is easy to see that
\[
\overline{Q}_1=Q_1,
\]
which implies that $Q_1$ is a real number.
Thus,
\[
\begin{array}{lll}
&&
-{a}_{11}(t)\big[ x_{12}(t) \overline{x}_{22}(t) x_{21}(t) \overline{x}_{11}(t)  + x_{11}(t)\overline{x}_{21} (t)x_{22}(t)  \overline{x}_{12} (t)\big]
\\
&&
-\big[ x_{12}(t) \overline{x}_{22}(t) x_{21}(t) \overline{x}_{11}(t)  + x_{11}(t)\overline{x}_{21} (t)x_{22}(t)  \overline{x}_{12} (t)\big]\overline{a}_{11}(t)
\\
&=&
-\big[ {a}_{11}(t)+\overline{a}_{11}(t)\big] \big[ x_{12}(t) \overline{x}_{22}(t) x_{21}(t) \overline{x}_{11}(t)  + x_{11}(t)\overline{x}_{21} (t)x_{22}(t)  \overline{x}_{12} (t)\big].
 \end{array}
 \]
It is obvious that
\[
 |x_{11}(t)|^2  a_{21}(t) x_{12}(t)\overline{x}_{22}(t) +|x_{11}(t)|^2  {x}_{22}(t) \overline {x}_{12}(t) \overline{a}_{21}(t)
 =2|x_{11}(t)|^2 \Re \{a_{21}(t) x_{12}(t)\overline{x}_{22}(t)\},
 \]
and
\[
|x_{11}(t)|^2  \overline{a}_{21}(t) {x}_{22}(t) \overline {x}_{12}(t) + x_{12}(t)\overline{x}_{22}(t)a_{21}(t) |x_{11}(t)|^2
=2|x_{11}(t)|^2 \Re \{\overline {a}_{21}(t) {x}_{22}(t) \overline{x}_{12}(t)\}.
\]
Now we claim that
\begin{equation}
\Re \{a_{21}(t) x_{12}(t)\overline{x}_{22}(t)\}=\Re \{\overline {a}_{21}(t) \overline{ x_{12}(t)\overline{x}_{22}(t)}\}.
\label{ReRe}
\end{equation}
In fact, for any $a,\, b\in \mathbb{H}$, it is easy to check that
\[
 \Re\{ a  \overline{b} \}= \Re \{ \overline{a} b\}.
 \]
Thus, it follows from (\ref{ReRe}) that
\[
\begin{array}{lll}
&&
 |x_{11}(t)|^2  a_{21}(t) x_{12}(t)\overline{x}_{22}(t) +|x_{11}(t)|^2  {x}_{22}(t) \overline {x}_{12}(t) \overline{a}_{21}(t)
 \\
 &&
 -|x_{11}(t)|^2  \overline{a}_{21}(t) {x}_{22}(t) \overline {x}_{12}(t) - x_{12}(t)\overline{x}_{22}(t)a_{21}(t) |x_{11}(t)|^2=0
 \end{array}
 \]
Similarly, we have
\[
\begin{array}{lll}
&&
 |x_{12}(t)|^2  a_{21}(t) x_{11}(t)\overline{x}_{21}(t) +|x_{12}(t)|^2  {x}_{21}(t) \overline {x}_{11}(t) \overline{a}_{21}(t)
 \\
 &&
 -|x_{12}(t)|^2  \overline{a}_{21}(t) {x}_{21}(t) \overline {x}_{11}(t) - a_{21}(t) x_{11}(t)\overline{x}_{21}(t) |x_{12}(t)|^2=0.
 \end{array}
 \]
Substituting above equalities into (\ref{Deri}), we have
\begin{equation}
\begin{array}{lll}
\displaystyle\frac{d}{dt}  W_{QDE}(t)
&=&
  \displaystyle \Big\{\big[ {a}_{11}(t) + {a}_{22}(t)  +\overline{a}_{11}(t)  + \overline{a}_{22}(t)  \big] |x_{11}(t)|^2|x_{22}(t)|^2
     \\
     &&
     +
        \big[ {a}_{11}(t) + {a}_{22}(t)  +\overline{a}_{11}(t)  + \overline{a}_{22}(t)  \big] |x_{12}(t)|^2|x_{21}(t)|^2
\\
&&
-  \big[ {a}_{11}(t) + {a}_{22}(t)  +\overline{a}_{11}(t)  + \overline{a}_{22}(t)  \big] x_{12}(t) \overline{x}_{22}(t) x_{21}(t) \overline{x}_{11}(t)  \\
&&
 - \big[ {a}_{11}(t) + {a}_{22}(t)  +\overline{a}_{11}(t)  + \overline{a}_{22}(t)  \big] x_{11}(t)\overline{x}_{21} (t)x_{22}(t)  \overline{x}_{12} (t)\big]
 \Big\}
 \\
 &=&
 \displaystyle  \big[ {a}_{11}(t) + {a}_{22}(t)  +\overline{a}_{11}(t)  + \overline{a}_{22}(t)  \big]   W_{QDE}(t)
   \\
 &=&
  \displaystyle \big[ {\rm tr} A(t) + {\rm tr} A^{+}(t)  \big]   W_{QDE}(t)
 .
\end{array}
\label{2.9}
\end{equation}
Integrating (\ref{2.9}) over $[t_0,t]$,  Liouville formula (\ref{Liou}) follows.
Consequently,  if $W_{QDE}(t)=0$ at some $t_0$ in $I$ then $W_{QDE}(t)=0$ on $I$.

\begin{remark}\label{Rem4.3}
It should be noted that the {\em Wronskian} of QDEs can also be defined by
\[
 W_{QDE}(t) =\rm{ddet} M(t): =  {\rm {rdet}} \big( M^+(t) M(t)\big).
 \label{W}
 \]
 By this definition, Proposition \ref{Pro4.1} and Theorem \ref{Th4.1} can be similarly proved.
\end{remark}

We need a lemma from (Theorem 2.10 \cite{Ky-AMC}).

\begin{lemma}\label{Lem4.1}
 The necessary and sufficient condition of invertibility of quaternionic matrix $M$ is $\rm{ddet} M\neq0$ (or $W_{QDE}(t)\neq 0$).
 \end{lemma}

\begin{proposition}\label{Pro4.2}
  If $W_{QDE}(t)=0$ at some $t_0$ in $I$ then  $x_1(t)$ and $x_2(t)$ are right linearly dependent on $I$.
  \end{proposition}

 {\bf PROOF.} From Liouville formula (Theorem \ref{Th4.1}), we have
 \[
 W_{QDE}(t_0)=0,\,\,\,\,\, \mbox{implies} \,\,\,\,\, W_{QDE}(t)=0,\,\,\,\, \mbox{for any}\,\,\,\,t\in I,
 \]
which implies that the quaternionic matrix $M$ is not invertible on $I$, in view of Lemma 4.1. Hence, the
 linear system
 \[
 M q=0,\,\,\,\, or \,\,\,\, (x_1(t),x_2(t))q=0,     \,\,\,\, q=(q_1,q_2)^T\in \mathbb{H},
 \]
 has a non-zero solution. Consequently, the two solution $x_1(t)$ and $x_2(t)$ are right linearly dependent on $I$.

 From Proposition \ref{Pro4.1}, Proposition \ref{Pro4.2} and Theorem \ref{Th4.1}, we immediately have

 \begin{theorem}\label{Th4.2}
 Let $A(t)$ in Eq. (\ref{2.1}) be continuous functions of $t$ on an interval $I$. Two solutions $x_1(t),x_2(t)$ of Eq. (\ref{2.1}) on $I$ are right linearly dependent on $I$ if and only if the absolute value of the Wronskian, $W_{QDE}(t)$, is zero at some $t_0$ in $I$.
 \end{theorem}

Now, we present an important result on the superposition principle and structure of the general solution.

\begin{theorem}\label{Th4.3} (Structure of the general solution) If  $x_1(t)$ and $x_2(t)$ are two solutions of Eq. (\ref{2.1}), then
\[
x(t)=  x_1(t) q_1 + x_2(t) q_2,
\]
is also a solution of (\ref{2.1}), where $q_1, \, q_2$ are quaternionic constants.
Let ${\cal {S}}_{\mathbb{H}}$ denote the set which contains of all the solutions of (\ref{2.1}). Then the set ${\cal {S}}_{\mathbb{H}}$ is a {\bf right
 free module} (the quaternionic vector space $\mathbb{H}^2$ over a division ring $\mathbb{H}$).
\end{theorem}

\begin{remark}\label{Rem-1D}
We consider the one-dimensional nonautonomous system
\[
\dot{x}=a(t) x.
\]
We see that it has a fundamental solution denoted by $\bar{x}(t)$. The set of solutions ${\cal {S}}=\{\bar{x}(t) q| q\in \mathbb{H}\}$ is a 1-dimensional right free module. In fact, $q \bar{x}(t)$ does not satisfy the equation, that is $q \bar{x}(t)\not\in {\cal {S}}$. Simple computation shows:
\[
\frac{d}{dt} (q \bar{x}(t))=q \dot\bar{x}(t)= q a(t) \bar{x}(t)\neq a(t) q\bar{x}(t).
\]
Consequently, $q \bar{x}(t)$ {\bf is not a solution of $\dot{x}=a(t) x$}, i.e. $q \bar{x}(t)\not\in {\cal {S}}$.

In particular, we consider the one-dimensional autonomous system
\[
\dot{x}=a  x.
\]
We see that it has a fundamental solution denoted by $\bar{x}(t)=\exp\{a t\}$. The set of solutions ${\cal {S}}=\{\exp\{a t\} q| q\in \mathbb{H}\} $ is a 1-dimensional right free module.

\end{remark}

Suppose that $x_1(t),\,x_2(t)$ are two fundamental solutions (they are independent in classical sense, i.e., the quaternionic vector space $\mathbb{H}^2$ over a field $\mathbb{C}$) of (\ref{2.1}). It is easy to see that the solution set $ {\cal{S}_\mathbb{C}}=\{x(t) \,\,\,\mbox{is a solution of (\ref{2.1})}\,\,\, \big|x(t)=c_1 x_1(t)+c_2x_2(t),\,c_1,c_2\in\mathbb{C}\}$ (the quaternionic vector space $\mathbb{H}^2$ over a field $\mathbb{C}$) is a linear vector space. However,
we claim that the set ${\cal{S}_\mathbb{C}}$ is not complete. The following counterexample shows this fact.\\
\noindent {\bf Couterexample:} Consider the following QDEs
\begin{equation}
\dot {x}= \left ( \begin{array}{ll}
\mathbf{i}  &  0
 \\
 0  &  \mathbf{j}
\end{array}
\right )x,\,\,\,\, x=(x_1,x_2)^T.
\label{CounterEX}
\end{equation}
It is easy to see that $x_1(t)=(e^{\mathbf{i}t},0)^{T}$, $x_2(t)=(0,e^{\mathbf{j}t})^{T}$ are two fundamental solutions of (\ref{CounterEX}). Also, it is easy to verify that $x(t)=(e^{\mathbf{i}t}\mathbf{k},0)^{T}$ is a solution of (\ref{CounterEX}).
Considering the solution set $ {\cal{S}_\mathbb{C}}=\{x(t) \,\,\,\mbox{is a solution of (\ref{2.1})}\,\,\, \big|x(t)=c_1 x_1(t)+c_2x_2(t),\,c_1,c_2\in\mathbb{C}\}$, we claim that the set $ {\cal{S}_\mathbb{C}}$ is not complete. If this is not the case, there exists the complex numbers $a,b$ such that $x(t)=x_1(t) a+x_2(t) b$, i.e., $(e^{\mathbf{i}t}\mathbf{k},0)^{T}=(e^{\mathbf{i}t},0)^{T} a+ (0,e^{\mathbf{j}t})^{T}b$. Then,
\begin{equation}
 \left\{ \begin{array}{ccc}
 e^{\mathbf{i}t}\mathbf{k} & =& e^{\mathbf{i}t}a,
 \\
0 &=&e^{\mathbf{j}t}b.
\end{array}
\right.
\label{EV1.1111}
\end{equation}
Consequently, we see that $-a\mathbf{k}=1,b=0$. Therefore, $a=\mathbf{k},b=0$, which contradicts to that $a$ is a complex number.

{\bf PROOF of THEOREM \ref{Th4.3}.}  The first assertion is easily obtained by direct differentiation.
Now we prove the second assertion. We proceed this with two steps. First,
we will prove the existence of two linearly independent solutions of Eq. (\ref{2.1}).
We can easily choose two linearly independent quaternionic vector $x_1^0=(x_{11}^0, x_{21}^0)^T\in \mathbb{H}^2$ and $x_2^0=(x_{12}^0, x_{22}^0)^T\in \mathbb{H}^2$ (e.g. the natural basis $e_1=(1, 0)^T, e_2=(0, 1)^T$).  According to Theorem 3.1, for any $k=1,2$ and $t_0\in I$, Eq. (\ref{2.1}) exists a unique solution $x_k(t),\,k=1,2$ satisfying $x_k(t_0)=x_k^0$.
If there are two quaternionic constants $q_1,q_2\in \mathbb{H}$ such that
\[
x_1(t) q_1 + x_2(t) q_2 =0,\,\,\,\,for\,\,\, any \,\,\,\, t\in I.
\]
In particularly, taking $t=t_0$, we have
\[
x_1(t_0) q_1 + x_2(t_0) q_2 =0,\,\,\,\,i.e.,\,\,\,\, x_1^0 q_1 + x_2^0 q_2 =0.
\]
Noting that $x_1^0$ and $x_2^0$ are two right linearly independent quaternionic vector, it follows that $q_1=q_2=0$. Thus,  $x_1(t),x_2(t)$ of Eq. (\ref{2.1}) on $I$ are right linearly dependent on $I$.

Secondly, we prove that any solution $x(t)$ of  Eq. (\ref{2.1}) can be represented by the  combination of above two right linearly independent solutions $x_1(t),x_2(t)$
\begin{equation}
x(t)=  x_1(t) q_1 + x_2(t) q_2,\,\,\,\,for \,\,\, any \,\,\,\, t\in I,
\label{Combination}
\end{equation}
where $q_1,q_2\in \mathbb{H}$.
On one hand, $x_1^0$ and $x_2^0$ are two right linearly independent quaternionic vector, which are a basis of two dimensional vector space $ \mathbb{H}^2$. Thus, there are two quaternionic constants $q_1,q_2\in \mathbb{H}$ such that
\[
x(t_0)=  x_1^0 q_1 + x_2^0 q_2 = x_1(t_0) q_1 + x_2(t_0) q_2.
\]
On the other hand, from the first assertion of this theorem, we know that
$
 x_1(t) q_1 + x_2(t) q_2
$
is also a solution of (\ref{2.1}) satisfying the IVP $x(t_0)$.  Therefore, $x(t)$ and $
 x_1(t) q_1 + x_2(t) q_2
$ are two solutions of  Eq. (\ref{2.1}) satisfying the same IVP $x(t_0)$. By the uniqueness theorem (Theorem \ref{Th3.1}), the equality (\ref{Combination}) holds.

\section{Fundamental matrix and solution to QDEs}

To study the solutions of Eq. (\ref{2.1}), it is necessary to
 introduce the quaternion-valued exponential function.
If $q\in \mathbb{H}$, the {\em exponential} of $q$ is defined by
\[
\exp q= \sum_{n=0}^{\infty} \frac{ q^n} {n!}
\]
(the series converges absolutely and uniformly on compact subsets).
If $r\in \mathbb{H}$ is
such that $qr = rq$, the exponential satisfies the addition theorem
\[
\exp\{ q\} \exp \{ r \} = \exp \{ q + r\}.
\]
\[
\exp\{ q\} = (\exp \Re q)[ (\cos \|\Im q\|) e + (\sin \| \Im q\|) \frac{\Im q} {\|\Im q\|}].
\]
Thus, $\exp q=1$ if and only if
\[
\Re q =0,\,\,\, and \,\,\, \|\Im q\|=0,\,\,\, ( mod 2\pi).
\]
Since $[q(t)r(t)]'=q'(t)r(t) + q(t) r'(t)$, ( here $'=\frac{ \rm d}{\rm dt}$),
\begin{equation}
[q^{n}(t)]'=\displaystyle\sum_{j=0}^{n-1} q^{j}(t) q'(t)q^{n-1-j}(t).
\label{nD}
\end{equation}

\begin{remark}\label{Rem5.1}
 It is known that if $q(t)$ reduces to the real-valued function,
\begin{equation}
[q^{n}(t)]'=nq^{n-1}(t),
\label{rnD}
\end{equation}
which is different from (\ref{nD}). This is the difference between quaternion-valued function and real function (complex function). But if $q$ is a quaternion function and $q(t)q'(t) = q'(t)q(t)$,  then  (\ref{rnD}) holds.
\end{remark}

\begin{proposition}\label{Pro5.1}
  If $q(t)$ is differentiable and if $q(t)q'(t) = q'(t)q(t)$, it
follows from that
\begin{equation}
[\exp q(t) ]' = [\exp q(t) ] q'(t).
\label{EXD}
\end{equation}
\end{proposition}

\begin{remark}\label{Rem5.2}
 If $q(t)\equiv q$ is a constant quaternion, then $q(t)q'(t) = q'(t)q(t)$ is always satisfied. We see the difference of the derivative between the quaternionic exponential function and complex (real) exponential function. For quaternion-value functions, the condition $q(t)q'(t) = q'(t)q(t)$ is needed to guarantee that the equality (\ref{EXD}) hold. For the real or complex function, $q(t)q'(t) = q'(t)q(t)$ is always satisfied.
\end{remark}

We prove Proposition \ref{Pro5.1} by two methods.

{\bf PROOF 1}: The way of Taylor expansion. Since $[q(t)r(t)]'=q'(t)r(t) + q(t) r'(t)$, by using (\ref{nD}),
\[
 \begin{array}{lll}
[\exp q(t) ]'&=&[\sum_{n=0}^{\infty} \frac{ q^n(t)} {n!}]'
=\displaystyle\sum_{n=0}^{\infty}\sum_{j=0}^{n-1} \frac{q^{j}(t) q'(t)q^{n-1-j}(t)}{n!}
\\
&=&
\displaystyle\sum_{n=0}^{\infty}\frac{   nq^{n-1}(t)} {n!}q'(t)
=\displaystyle\sum_{n=0}^{\infty}\frac{    q^{n-1}(t)} {(n-1)!}q'(t)=[\exp q(t) ] q'(t).
 \end{array}
\]

{\bf PROOF 2}: The way of Euler formula. We rewrite $q(t)$ as
\[
q(t)=q_0(t) + q_1(t) \mathbf{i} + q_2 (t) \mathbf{j} + q_3 (t) \mathbf{k} = q_0(t) + \Im {q}(t),
\]
where
\[
\Im {q}(t)=q_1(t) \mathbf{i} + q_2 (t) \mathbf{j} + q_3 (t) \mathbf{k} = \frac{\Im {q}(t)}{|\Im {q}(t)|}|\Im {q}(t)|.
\]
Then we can rewrite $\exp \{ q(t)\}$ as
\[
\begin{array}{lll}
 e^{q(t)}
 &=&
  e^{q_0(t)+\Im {q}(t)}
  \\
  &=&
  e^{q_0(t)}e^{\Im {q}(t)}
    \\
  &=&
  e^{q_0(t)}e^{\frac{\Im {q}(t)}{|\Im {q}(t)|}|\Im {q}(t)|}
      \\
  &=&
  e^{q_0(t)}[\cos (|\Im {q}(t)|) + \frac{\Im {q}(t)}{|\Im {q}(t)|}\sin( |\Im {q}(t)|)]
  \end{array}
  \]
If $q(t)=0$, obviously, the equality (\ref{EXD}) holds.
If $q(t)\neq 0$, differentiate $\exp \{ q(t)\}$ with respect to $t$, we have
\[
\begin{array}{lll}
&&
\displaystyle\frac{\rm d}{\rm dt} e^{q(t)}
\\
&=&
\displaystyle\frac{\rm d}{\rm dt}\Big\{ e^{q_0(t)}[\cos (|\Im {q}(t)|) + \frac{\Im {q}(t)}{|\Im {q}(t)|}\sin( |\Im {q}(t)|)]\Big\}
\\
&=&
\displaystyle q_0'(t)e^{q_0(t)}[\cos (|\Im {q}(t)|) + \frac{\Im {q}(t)}{|\Im {q}(t)|}\sin( |\Im {q}(t)|)]
\\
&&
+
\displaystyle e^{q_0(t)}\frac{\rm d}{\rm dt}[ \cos (|\Im {q}(t)|) + \frac{\Im {q}(t)}{|\Im {q}(t)|}\sin( |\Im {q}(t)|)]
\\
&=&
\displaystyle q_0'(t)e^{q_0(t)} \Big [\cos (|\Im {q}(t)|) + \frac{\Im {q}(t)}{|\Im {q}(t)|}\sin( |\Im {q}(t)|) \Big]
\\
&&
+
\displaystyle e^{q_0(t)}\Big[- \sin (|\Im {q}(t)|) \frac{\rm d}{\rm dt} |\Im {q}(t)|
 + \frac{\Im {q}'(t)|\Im {q}(t)|- \Im {q}(t) \frac{\rm d}{\rm dt} |\Im {q}(t)|}{|\Im {q}(t)|^2}\sin( |\Im {q}(t)|)
\\
&&
+
\hspace{1cm}
\displaystyle \frac{\Im {q}(t)}{|\Im {q}(t)|}\cos( |\Im {q}(t)|) \frac{\rm d}{\rm dt} |\Im {q}(t)|
 \Big]
\end{array}
\]
Since $q(t)q'(t) = q'(t)q(t)$, $ \frac{\Im {q}'(t)|\Im {q}(t)|- \Im {q}(t) \frac{\rm d}{\rm dt} |\Im {q}(t)|}{|\Im {q}(t)|^2}=0$, it follows that
\[
\begin{array}{lll}
&&
\displaystyle\frac{\rm d}{\rm dt} e^{q(t)}
\\
&=&
\displaystyle q_0'(t)e^{q_0(t)} \Big [\cos (|\Im {q}(t)|) + \frac{\Im {q}(t)}{|\Im {q}(t)|}\sin( |\Im {q}(t)|) \Big ]
\\
&&
+
\displaystyle e^{q_0(t)}\Big[- \sin (|\Im {q}(t)|) \frac{\rm d}{\rm dt} |\Im {q}(t)|
+ \frac{\Im {q}(t)}{|\Im {q}(t)|}\cos( |\Im {q}(t)|) \frac{\rm d}{\rm dt} |\Im {q}(t)|
 \Big]
 \\
&=&
\displaystyle q_0'(t)e^{q_0(t)} \Big [\cos (|\Im {q}(t)|) + \frac{\Im {q}(t)}{|\Im {q}(t)|}\sin( |\Im {q}(t)|) \Big]
\\
&&
+
\displaystyle e^{q_0(t)}\frac{\Im {q}(t)}{|\Im {q}(t)|}  \frac{\rm d}{\rm dt} |\Im {q}(t)|
\Big[-  \Im {q}^{-1}(t) |\Im {q}(t)| \sin (|\Im {q}(t)|)
+ \cos( |\Im {q}(t)|)
 \Big].
\end{array}
\]
Note that
\[
\Im {q}^{-1}(t) =\frac{\overline \Im {q}(t)}{|\Im {q}(t)|^2},
\]
and $-\overline \Im {q}(t)=\Im {q}(t)$. Thus, we have
\begin{equation}
\begin{array}{lll}
&&
\displaystyle\frac{\rm d}{\rm dt} e^{q(t)}
 \\
&=&
\displaystyle q_0'(t)e^{q_0(t)} \Big [\cos (|\Im {q}(t)|) + \frac{\Im {q}(t)}{|\Im {q}(t)|}\sin( |\Im {q}(t)|) \Big]
\\
&&
+
\displaystyle \frac{\rm d}{\rm dt} |\Im {q}(t)|\frac{\Im {q}(t)}{|\Im {q}(t)|}  e^{q_0(t)}
\Big[
\cos( |\Im {q}(t)|)+ \frac{\Im {q}(t)}{|\Im {q}(t)|}  \sin (|\Im {q}(t)|)
 \Big].
\end{array}
\label{De}
\end{equation}
From $q(t)q'(t) = q'(t)q(t)$, we conclude that $\Im q(t) \Im q'(t) = \Im q'(t) \Im q(t)$. Substituting $\Im {q}(t)$ by $q_1(t) i + q_2 (t) j + q_3 (t) k$, we obtain
\[
\begin{array}{lll}
q_1(t)q_2'(t) &=& q_1'(t)q_2(t),
\\
q_1(t)q_3'(t) &=& q_1'(t)q_3(t),
\\
q_2(t)q_3'(t) &=& q_2'(t)q_3(t).
\end{array}
\]
By these equalities, we have
\[
\begin{array}{lll}
\displaystyle \frac{\rm d}{\rm dt} |\Im {q}(t)|\frac{\Im {q}(t)}{|\Im {q}(t)|}
&=&
\displaystyle \frac{q_1(t)q_1'(t)+q_2(t)q_2'(t)+q_3(t)q_3'(t)}{\sqrt{q_1^2(t)+q_2^2(t)+q_2^2(t)}}\displaystyle \frac{q_1(t) \mathbf{i}+q_2(t)\mathbf{j}+q_3(t)\mathbf{k}}{\sqrt{q_1^2(t)+q_2^2(t)+q_2^2(t)}}
\\
&=&
\displaystyle \frac{[q_1^2(t)+q_2^2(t)+q_2^2(t)][q_1'(t) \mathbf{i}+q_2'(t)\mathbf{j}+q_3'(t)\mathbf{k}]}{  q_1^2(t)+q_2^2(t)+q_2^2(t) }
\\
&=&
\Im q'(t).
\end{array}
\]
Consequently, it follows from (\ref{De}) that
\[
\begin{array}{lll}
&&
\displaystyle\frac{\rm d}{\rm dt} e^{q(t)}
 \\
&=&
\displaystyle q_0'(t)e^{q_0(t)} \Big [\cos (|\Im {q}(t)|) + \frac{\Im {q}(t)}{|\Im {q}(t)|}\sin( |\Im {q}(t)|) \Big]
\\
&&
+
\displaystyle \Im q'(t)  e^{q_0(t)}
\Big[
\cos( |\Im {q}(t)|)+ \frac{\Im {q}(t)}{|\Im {q}(t)|}  \sin (|\Im {q}(t)|)
 \Big]
 \\
 &=&(q_0(t)+\Im q'(t))e^{q(t)}
 \\
 &=&q'(t)e^{q(t)}.
\end{array}
\]
This completes the proof of Proposition \ref{Pro5.1}.

If a matrix $A\in \mathbb{H}^{2\times 2}$, the {\em exponential} of $A$ is defined by
\[
\exp A= \sum_{n=0}^{\infty} \frac{ A^n} {n!}=E+\frac{ A} {1!}+\frac{ A^2} {2!}+\cdots.
\]
Moreover, for $t\in \mathbb{R}$,
\[
\exp (A t)= \sum_{n=0}^{\infty} \frac{ A^nt^n} {n!}=E+\frac{ A t} {1!}+\frac{ A^2 t^2} {2!}+\cdots.
\]

\begin{proposition}\label{Pro5.2}
 The sires $\exp A= \sum_{n=0}^{\infty} \frac{ A^n} {n!}$ is absolutely convergent. The sires $\exp (A t)= \sum_{n=0}^{\infty} \frac{ A^n} {n!}$ is uniformly convergent on any finite interval $I$.
 \end{proposition}

{\bf PROOF.}  First, we show that the sires $\exp A= \sum_{n=0}^{\infty} \frac{ A^n} {n!}$ is absolutely convergent. In fact, for any positive integer $k$, we have
\[
\|\frac{A^k}{k!}\|\leq \frac{\|A\|^k}{k!}.
\]
For any quaternionic matrix $A$, $\|A\|$ is a real number. Then it is easy to see that the real series $\sum_{n=0}^{\infty} \frac{ \|A\|^n} {n!}$ is convergent. Thus, the sires $\exp A= \sum_{n=0}^{\infty} \frac{ A^n} {n!}$ is absolutely convergent.

Now we prove the second conclusion. In fact, for any finite interval $I$, without loss of generality, we assume that $t\leq \mu$. Then we have
 \[
\|\frac{(At)^k}{k!}\|\leq \frac{\|A\|^k\|t\|^k}{k!}\leq \frac{\|A\|^k\|\mu\|^k}{k!}.
\]
 Note that $\frac{\|A\|^k\|\mu\|^k}{k!}$ is convergent. Thus, the sires $ \sum_{n=0}^{\infty} \frac{ A^n} {n!}$ is uniformly convergent on any finite interval $I$.

\begin{proposition}\label{Pro5.3}
 For two matrix $A,B\in \mathbb{H}^{2\times2}$, if $AB=BA$, then
\[
\exp(A+B)=\exp A\exp B.
\]
\end{proposition}

From Proposition \ref{Pro5.1}, we have

\begin{proposition}\label{Pro5.4}
If $A(t)$ is differentiable and if $A(t)A'(t) = A'(t)A(t)$, then
\[
[\exp A(t) ]' = [\exp A(t) ] A'(t).
\]
\end{proposition}

\begin{definition}\label{Def5.1}
  Let $x_1(t),x_2(t)$ be any two solutions of Eq. (\ref{2.1}) on $I$. Then we call
 \[
M(t)=(x_1(t),x_2(t))^T
=
\left ( \begin{array}{ll}
  {x} _{11}(t)  & {x} _{12}(t)
 \\
 {x} _{21}(t)  & {x} _{22}(t)
\end{array}
\right )
\]
as a {\em solution matrix} of Eq. (\ref{2.1}). Moreover, if $x_1(t),x_2(t)$ be two right linearly independent solutions of Eq. (\ref{2.1}) on $I$, the solution matrix of Eq. (\ref{2.1}) is said to be a {\em fundamental matrix} of Eq. (\ref{2.1}). Further, if $M(t_0)=E(identity)$, $M(t)$ is said to be a {\em normal fundamental matrix}.
\end{definition}

\begin{remark}\label{Rem5.3}
  Note that fundamental matrix is not unique. From Definition \ref{Def5.1}, it is easy to know that if $M(t)$ is a {\em solution matrix} of Eq. (\ref{2.1}). $M(t)$ also satisfies Eq. (\ref{2.1}), that is,
\[
\dot{M}(t)=A(t) M(t).
\]
\end{remark}

From Theorem \ref{Th4.2}, we know that if $M(t)$ is a fundamental matrix of Eq. (\ref{2.1}), the Wronskian determinant $W_{QDE}(t)\neq0$. By Theorem \ref{Th4.3}, we have

\begin{theorem}\label{Th5.1}
  Let $M(t)$ be a fundamental matrix of Eq. (\ref{2.1}). Any solution $x(t)$ of  Eq. (\ref{2.1}) can be represented by
\[
x(t)=M(t) q,
\]
where $q$ is a constant quaternionic vector. Moreover, for given IVP $x(t_0)=x^0$,
\[
x(t)=M(t)M^{-1}(t_0)x^0.
\]
\end{theorem}

By Theorem \ref{Th4.2}, we also have

\begin{theorem}\label{Th5.2}
 A solution matrix $M(t)$ of Eq. (\ref{2.1}) on $I$ is a fundamental matrix if and only if $\rm{ddet} M(t)\neq0$ (or $W_{QDE}(t)\neq 0$) on $I$. Moreover, for some $t_0$ in $I$ such that $\rm{ddet} M(t_0)\neq0$ (or $W_{QDE}(t_0)\neq0$ ), then $\rm{ddet} M(t)\neq0$ (or $W_{QDE}(t)\neq0$ ).
 \end{theorem}

\noindent {\bf Example 1} Show that
\[
M(t)=(\phi_1(t),\phi_2(t))^T
=
\left ( \begin{array}{ll}
e^{\boldsymbol{i} t}  &  e^{\boldsymbol{i} t}\int_{t_{0}}^{t}e^{-\boldsymbol{i} t}\cdot e^{\boldsymbol{j} t}dt
\\
0  & \,\,\,\,\,\,\,\,\,\,\,\,\,\, e^{\boldsymbol{j} t}
\end{array}
\right )
\]
is a fundamental matrix of the QDES
\begin{equation}
 \left ( \begin{array}{ll}
 \dot {x}_1
 \\
\dot {x}_2
\end{array}
\right )
=
\left ( \begin{array}{ll}
 \mathbf{i}  &  1
 \\
 0  &  \mathbf{j}
\end{array}
\right ) \left ( \begin{array}{ll}
   {x}_1
 \\
  {x}_2
\end{array}
\right ) .
\label{Ex5.1}
\end{equation}

{\bf Proof.} First, we show that $M(t)$ is a solution matrix of Eq.(\ref{Ex5.1}). In fact, let $\phi_1(t)=\left ( \begin{array}{ll}
  e^{\boldsymbol{i} t}
 \\
 0
\end{array}
\right )$, in view of $ (\boldsymbol{i} t)' (\boldsymbol{i}t)= (\boldsymbol{i}t) (\boldsymbol{i}t )'$, then
\[
\dot{\phi}_1(t)=\left ( \begin{array}{ll}
\boldsymbol{i}  e^{\boldsymbol{i} t}
 \\
 0
\end{array}
\right )
=
\left ( \begin{array}{ll}
 \boldsymbol{i}  &  1
 \\
 0  &  \boldsymbol{j}
\end{array}
\right )
\left ( \begin{array}{ll}
e^{\boldsymbol{i} t}
 \\
 0
\end{array}
\right )
=
\left ( \begin{array}{ll}
 \boldsymbol{i}  &  1
 \\
 0  &  \boldsymbol{j}
\end{array}
\right ) \phi_1(t),
\]
which implies that $\phi_1(t)$ is a solution of Eq.(\ref{Ex5.1}).
Similarly, let $\phi_2(t)=\left ( \begin{array}{ll}
 e^{\boldsymbol{i} t}\int_{t_{0}}^{t}e^{-\boldsymbol{i} t}\cdot e^{\boldsymbol{j} t}dt
 \\
\,\,\,\,\,\,\,\,\,\,\,\,\,\, e^{\boldsymbol{j} t}
\end{array}
\right )$, then
\[
\dot{\phi}_2(t)=\left ( \begin{array}{ll}
\boldsymbol{i}e^{\boldsymbol{i} t}\int_{t_{0}}^{t}e^{-\boldsymbol{i} t}\cdot e^{\boldsymbol{j} t}+e^{\boldsymbol{j} t}dt
 \\
\,\,\,\,\,\,\,\,\,\,\,\,\,\,\,\,\,\,\,\,\,\,  \boldsymbol{j} e^{\boldsymbol{j} t}
\end{array}
\right )
=
\left ( \begin{array}{ll}
 \boldsymbol{i}  &  1
 \\
 0  &  \boldsymbol{j}
\end{array}
\right )
\left ( \begin{array}{ll}
 e^{\boldsymbol{i} t}\int_{t_{0}}^{t}e^{-\boldsymbol{i} t}\cdot e^{\boldsymbol{j} t}dt
 \\
\,\,\,\,\,\,\,\,\,\,\,\,\,\,e^{\boldsymbol{j} t}
\end{array}
\right )
=
\left ( \begin{array}{ll}
 \boldsymbol{i}  &  1
 \\
 0  &  \boldsymbol{j}
\end{array}
\right ) \phi_2(t),
\]
 which implies that $\phi_2(t)$ is another solution of Eq.(\ref{Ex5.1}). Therefore,  $M(t)=(\phi_1(t),\phi_2(t))^\top$ is a solution matrix of Eq.(\ref{Ex5.1}).

Secondly, by Theorem \ref{Th5.2}, and taking $t_0=0$,
\[
\rm{ddet} M(t_0)= \rm{rdet} (M(t_0) M^{+}(t_0)) \neq0,
\]
then $\rm{ddet} M(t)\neq0$. Therefore, $M(t)$ is a fundamental matrix of Eq. (\ref{Ex5.1}).

From Theorem \ref{Th5.1} and Theorem \ref{Th5.2}, the following corollaries follows immediately.

\begin{corollary}\label{Cor5.1}  Let $M(t)$ be a fundamental matrix of Eq. (\ref{Ex5.1}) on $I$. $Q\in \mathbb{H}^{2\times 2}$ is an invertible constant quaternionic matrix. Then $M(t)Q$ is also a fundamental matrix of Eq. (\ref{Ex5.1}) on $I$.
\end{corollary}

{\bf PROOF.} As pointed out in Remark \ref{Rem5.3}, any solution matrix $M(t)$ satisfies Eq. (\ref{2.1}), that is,
\[
\dot{M}(t)=A(t) M(t).
\]
The inverse also holds.
Taking $\Psi(t)=M(t)Q$, differentiating it, and noticing that $M(t)$ is a fundamental matrix satisfying Eq. (\ref{2.1}), we have
\[
\Psi'(t)= M'(t)Q= A(t)M(t)Q=A(t)\Psi(t),
\]
which implies that $\Psi(t)$ is a solution matrix of Eq. (\ref{2.1}). Moreover, since $Q$ is invertible, we have
\[
\rm{ddet} \Psi(t)=\rm{ddet} M(t)Q= \rm{ddet} M(t)\cdot \rm{ddet}Q \neq0.
\]
Therefore, $\Psi(t)=M(t)Q$ is also a fundamental matrix of Eq. (\ref{2.1}).

\begin{corollary}\label{Cor5.2}
 If $M(t),\Psi(t)$ are two fundamental matrices of Eq. (\ref{2.1}) on $I$,  there exists an invertible constant quaternionic matrix $Q\in \mathbb{H}^{2\times 2}$ such that
\[
\Psi(t)=M(t)Q,\,\,\,t \in I.
\]
\end{corollary}
\indent {\bf PROOF.}  Since $M(t)$ is a fundamental matrix, $M(t)$ is invertible and denote the inverse by $M^{-1}(t)$. Set
\begin{equation}
\Psi(t)=M(t) X(t),\,\,t\in I.
\label{X}
\end{equation}
$M(t),\Psi(t)$ are fundamental matrices, so ${\rm ddet}M(t)\neq0$ and  ${\rm ddet}\Psi(t)\neq0$. Thus, ${\rm ddet}X(t)\neq0$.
Differentiation on (\ref{X}) gives
\[
\Psi'(t)=M'(t)X(t)+M(t)X'(t)=A(t)M(t)X(t)+M(t)X'(t)=A(t)\Psi(t)+M(t)X'(t)
\]
On the other hand, $\Psi(t)$ is a fundamental matrix, which implies that
\[
\Psi'(t)=A(t)\Psi(t).
\]
Compare above two equalities, we have
\[
M(t)X'(t)=0.
\]
Noticing that $M(t)$ is invertible, we obtain that $X'(t)=0$, which implies that $X(t)$ is a constant matrix. Letting $X(t)=Q$. The proof is complete.

\smallskip

Now consider the following quaternion-valued equations with constant quaternionic coefficients
\begin{equation}
\dot{x}(t) = A x(t),
\label{Constant}
\end{equation}
where $A\in   \mathbb{H}^{2\times 2} $ is a  constant quaternion matrix. Then we have

\begin{theorem}\label{Th5.3}
  $M(t)= \exp\{A t\}$ is a fundamental matrix of (\ref{Constant}). Moreover, any solution $x(t)$ of (\ref{Constant}) can be represented by
\[
x(t)=\exp\{A t\} q,
\]
where $q\in   \mathbb{H}^{2} $ is an arbitrary constant quaternion. For the IVP $x(t_0)=x^0$, any solution $x(t)$ of (\ref{Constant}) can be represented by
\[
x(t)=\exp\{A( t-s)\} x^0.
\]
\end{theorem}

 {\bf Proof.} Note that $A\in   \mathbb{H}^{2\times 2} $ is a  constant quaternion matrix, by using
 $
A t \cdot (A t)'= (A t)' \cdot A t
$
we have
 \[
M'(t)= [\exp\{A t\}]'=A \exp\{A t\}=A M(t),
 \]
  which implies $ \exp\{A t\}$ is a solution matrix of (\ref{Constant}). Moreover, $\rm{ddet} M(0)=\rm{ddet} \exp(A\cdot 0)=\rm{ddet} E=1\neq 0$. Consequently,
 $ \rm{ddet} M(t)\neq 0$. Therefore, $ \exp\{A t\}$ is a fundamental matrix of (\ref{Constant}).

Now consider the diagonal homogenous system.
\begin{equation}
\left ( \begin{array}{ll}
 \dot{x} _{1}(t)
 \\
\dot{x} _{2}(t)
\end{array}
\right )
=
\left ( \begin{array}{ccc}
  {a} _{1}(t)  & 0
 \\
0 & {a} _{2}(t)
\end{array}
\right )
\left ( \begin{array}{ll}
  {x} _{1}(t)
 \\
 {x} _{2}(t)
\end{array}
\right )
\label{diag}
\end{equation}

\begin{theorem}\label{Th5.4}
  Assume that
\begin{equation}
a_i(t)   \int_{t_0}^{t} a_i(s) ds  =    \int_{t_0}^{t} a_i(s) ds  a_i(t)
\label{Condition1}
\end{equation}
holds. The fundamental matrix can be chose by
\[
M(t)=
\left ( \begin{array}{ccc}
  \exp\{ \int_{t_0}^{t} a_1(s) ds\}   & 0
 \\
0 &            \exp\{ \int_{t_0}^{t} a_2(s) ds\}
\end{array}
\right ).
\]
Then the solution of the diagonal system (\ref{diag}) with the initial value $x(t_0)=x^0$ is given by
\[
x(t)=M(t)x^0.
\]
\end{theorem}

\section{\bf Algorithm for computing fundamental matrix 
}

Since the fundamental matrix plays great role in solving QDEs, in this section, we provide an algorithm for computing fundamental matrix of linear QDEs with constant coefficients.

\subsection{ Method 1: using expansion of $\exp \{At\}$}

From Theorem \ref{Th5.3}, we know that $\exp\{A t\}$ is a fundamental matrix of (\ref{Constant}).
So if the coefficient matrix is not very complicate, we can use the definition of $\exp (At)$ to compute fundamental matrix of linear QDEs with constant coefficients.

\begin{theorem}\label{Th6.1}
 If $A=diag( \lambda_1,\lambda_2)\in \mathbb{H}^{2\times2}$ is a diagonal matrix, then
\[
\exp \{At \}= \left ( \begin{array}{ccc}
  \exp\{  \lambda_1t\}    &0
  \\

     0&       \exp\{ \lambda_2 t\}
\end{array}
\right ).
\]
\end{theorem}

{\bf Proof.} By the expansion,
\[
\exp \{ At \}=E+  \left ( \begin{array}{ccc}
   \lambda_1     &0
  \\

     0&       \lambda_2
\end{array}
\right )\frac{t}{1!}
+
 \left ( \begin{array}{ccc}
   \lambda_1     &0
  \\

     0&       \lambda_2
\end{array}
\right )^2\frac{t^2}{2!}
+\cdots
=
\left ( \begin{array}{ccc}
  \exp\{  \lambda_1t\}    &0
  \\

     0&       \exp\{ \lambda_2 t\}
\end{array}
\right ).
\]

By Proposition \ref{Pro5.3}, we can divide the matrix to some simple ones and use the expansion to compute the fundamental matrix.
\[
A=diag A+ N,
\]
where $N$ is a nilpotent matrix. That is, $N^n=0$ and $n$ is a finite number.

\smallskip

\noindent {\bf Example 2} Find a fundamental matrix of the following QDES
\[
\dot {x}= \left ( \begin{array}{ll}
 \mathbf{k}  &  1
 \\
 0  &  \mathbf{k}
\end{array}
\right )x,\,\,\,\, x=(x_1,x_2)^T.
\]

{\bf Answer.} We see that
\[
\left ( \begin{array}{ll}
 \mathbf{k}  &  1
 \\
 0  &  \mathbf{k}
\end{array}
\right )
=\left ( \begin{array}{ll}
 \mathbf{k}  &  0
 \\
 0  &  \mathbf{k}
\end{array}
\right )+\left ( \begin{array}{ll}
 0  &  1
 \\
 0  &  0
\end{array}
\right ).
\]
Noticing that $\left ( \begin{array}{ll}
 \mathbf{k}  &  0
 \\
 0  &  \mathbf{k}
\end{array}
\right )
\left ( \begin{array}{ll}
 0  &  1
 \\
 0  &  0
\end{array}
\right )
=
\left ( \begin{array}{ll}
 0  &  1
 \\
 0  &  0
\end{array}
\right )\left ( \begin{array}{ll}
 \mathbf{k}  &  0
 \\
 0  &  \mathbf{k}
\end{array}
\right )
$, by Proposition 5.3 and Theorem 6.1, we have
\[
\begin{array}{lll}
\exp(At)
&=&
\exp\left ( \begin{array}{ll}
 \mathbf{k}  &  0
 \\
 0  &  \mathbf{k}
\end{array}
\right )
\exp\left ( \begin{array}{ll}
 0  &  1
 \\
 0  &  0
\end{array}
\right )
\\
&=&
\left ( \begin{array}{ll}
 e^{\mathbf{k}t}  &  0
 \\
 0  &  e^{\mathbf{k}t}
\end{array}
\right )
\Big(
E+  \left ( \begin{array}{ccc}
0     & 1
  \\
 0  &   0
\end{array}
\right )\frac{t}{1!}
+
 \left ( \begin{array}{ccc}
 0   &  1
  \\
 0  &   0
\end{array}
\right )^2\frac{t^2}{2!}
+\cdots
\Big)
.
\end{array}
\]
Note that
\[
 \left ( \begin{array}{ccc}
  0     & 1
  \\
  0   &   0
\end{array}
\right )^2
= \left ( \begin{array}{ccc}
  0     & 0
  \\
  0   &   0
\end{array}
\right ),\,\, \left ( \begin{array}{ccc}
  0     & 1
  \\
  0  &   0
\end{array}
\right )^3
= \left ( \begin{array}{ccc}
  0     & 0
  \\
     0&   0
\end{array}
\right ) ,\cdots,
\]
Then the fundamental matrix
\[
\exp \{ At \}
=
\left ( \begin{array}{ll}
 e^{\mathbf{k}t}  &  0
 \\
 0  &  e^{\mathbf{k}t}
\end{array}
\right ) \left ( \begin{array}{ll}
1  &  t
 \\
0  &  1
\end{array}
\right )
 =
\left ( \begin{array}{ll}
 e^{\mathbf{k}t}  &   t e^{\mathbf{k}t}
 \\
  0  &  e^{\mathbf{k}t}
\end{array}
\right ).
\]
For the matrix $\left ( \begin{array}{ll}
\lambda  &  1
 \\
 0  &  \lambda
\end{array}
\right )$, $\lambda\in\mathbb{H}$, similar computation shows that

 \begin{theorem}\label{Th6.1.1}
 For any $\lambda\in\mathbb{H}$, we have
\[
\exp \left\{ \left ( \begin{array}{ll}
\lambda  &  1
 \\
 0  &  \lambda
\end{array}
\right )t \right\}
=
\left ( \begin{array}{ll}
 e^{\lambda t}  &  0
 \\
 0  &  e^{\lambda t}
\end{array}
\right ) \left ( \begin{array}{ll}
1  &  t
 \\
0  &  1
\end{array}
\right ).
\]
\end{theorem}
If we generalize this result to $n$-dimensional case, we have the following formula.
 \begin{theorem}\label{Th6.1.2}
 For any $\lambda\in\mathbb{H}$, we have
\[
\exp \left\{ \left ( \begin{array}{ccccc}
\lambda & 1 &0 \cdots 0 & 0
\\
0 & \lambda &1 \cdots 0 & 0
\\
\,\,& \,\,& \,\,\,\,\,\,\, \ddots \,\,\ddots  &
\\
0 & 0 &0\cdots \lambda & 1
\\
0 & 0 &0\cdots 0& \lambda
\end{array}
\right )t \right\}
=
diag(
\exp{\lambda t}, \cdots, \exp{\lambda t}) \left ( \begin{array}{ccccc}
1 & t &\frac{t^{2}}{2!} \cdots \frac{t^{k-2}}{(k-2)!} & \frac{t^{k-1}}{(k-1)!}
\\
0 & 1 &t \cdots \frac{t^{k-3}}{(k-3)!} & \frac{t^{k-2}}{(k-2)!}
\\
\,\,& \,\,& \,\,\,\,\, \ddots \,\,\ddots  &
\\
0 & 0 &0 \cdots 1 & t
\\
0 & 0 &0 \cdots 0 & 1
\end{array}
\right ).
\]
\end{theorem}

We remark that this method can not be applied extensively when the coefficient matrix $A$ is complicate. Usually, the two divisions of $A$ can not be communicate. For example,
\[
\left ( \begin{array}{ll}
 \mathbf{i}  &  \mathbf{k}
 \\
 0  &  \mathbf{j}
\end{array}
\right )
=\left ( \begin{array}{ll}
 \mathbf{i}  &  0
 \\
 0  &  \mathbf{j}
\end{array}
\right )+\left ( \begin{array}{ll}
 0  &  \mathbf{k}
 \\
 0  &  0
\end{array}
\right ).
\]
We see that
\[
\left ( \begin{array}{ll}
 \mathbf{i}  &  0
 \\
 0  &  \mathbf{j}
\end{array}
\right )
\left ( \begin{array}{ll}
 0  &  \mathbf{k}
 \\
 0  &  0
\end{array}
\right )
\neq
\left ( \begin{array}{ll}
 0  &  \mathbf{k}
 \\
 0  &  0
\end{array}
\right )\left ( \begin{array}{ll}
 \mathbf{i}  &  0
 \\
 0  &  \mathbf{j}
\end{array}
\right )
\]
In this case, this method can not be used. Even, for a simple matrix $\left ( \begin{array}{ll}
 \mathbf{i}  &  \mathbf{1}
 \\
 0  &  \mathbf{j}
\end{array}
\right )$, you can not use this method. So we should find more effective method to compute the fundamental matrix. In what follows,
we try to use the eigenvalue and eigenvector theory to compute it.

\subsection {Method 2:  eigenvalue and eigenvector theory}

Since quaternion is a noncommunicative algebra, eigenvalue of a matrix $A\in \mathbb{H}^{2\times 2}$ should be defined by left eigenvalue and right eigenvalue, respectively. A quaternion $\lambda$ is said to be a right eigenvalue of $A$ if
\[
Ax=x \lambda
\]
 for
some nonzero (column) vector $x$ with quaternion components. Similarly, a
quaternion $\lambda$ is a left eigenvalue of $A$ if
\[
Ax= \lambda  x
\]
 for some nonzero (column) vector
$x$ with quaternion components. Right and left eigenvalues are in general unrelated. Usually, right and left eigenvalues are different. It should be noted that there are many differences in the eigenvalue problem for complex (real) and quaternion matrices. Eigenvalues of complex matrices
satisfy Brauer's theorem \cite{Br} for the inclusion of the eigenvalues. While right eigenvalues of quaternion matrices do not have this property. Another contrast in the eigenvalue problems for complex and
quaternion matrices is that a complex $n\times n$ matrix cannot have more than $n$ complex
eigenvalues, while it can have infinitely many quaternion left eigenvalues (see Theorem 2.1 in \cite{ZhangFZ}). Finally, it is well-known that any eigenvalue of an $n\times n$
complex matrix must be a root of its characteristic polynomial. Left eigenvalues
and right eigenvalues of quaternion matrices do not share this property of
complex matrices.

In this paper, we focus on the right eigenvalue and eigenvector. Because we emphasize on finding the solution taking the form
\[
x=q e^{\lambda t},\,\,\,q=(q_1,q_2)^T,
\]
where
$\lambda$ is a quaternionic constant and $q$ is a constant quaternionic vector. Substituting it into Eq. (\ref{Constant}), we have
\[
q \lambda  e^{\lambda t} = A q  e^{\lambda t}.
\]
Because $ e^{\lambda t}\neq 0$, it follows that
\[
q \lambda = A q,
\]
or
\begin{equation}
 A q=q \lambda.
\label{CQ}
\end{equation}
So if we find such eigenvalue $\lambda$ and eigenvector $q$ in (\ref{CQ}), we will find the solution of Eq. (\ref{Constant}). We also say that the eigenvalue $\lambda$ and eigenvector $q$ are the right eigenvalue and eigenvector of the QDEs (\ref{Constant}).

Notice also that if $0\neq \alpha\in \mathbb{H}$, then
\[
A q = q \lambda \Rightarrow  A q \alpha = q \alpha (\alpha^{-1} \lambda \alpha),
\]
so we can sensibly talk about the eigenline spanned by an eigenvector $q$, even though there may
be many associated eigenvalues! In what follows, if
\[
\theta=\alpha^{-1} \lambda \alpha
\]
we call it that $\theta$ is similar to $\lambda$. If there exists a nonsingular matrix $T$ such that
\[
A=T^{-1}BT,
\] we call it that $A$ is similar to $B$, denoted by $A\sim B$.

\begin{remark}\label{Rem6.1}
 If $\lambda$ is a characteristic root of $A$, then so is $\alpha^{-1} \lambda \alpha$.
 \end{remark}

From the definition of fundamental matrix, we have

\begin{theorem}\label{Th6.2}
 If the matrix $A$ has two right linearly independent eigenvectors $q_1$ and $q_1$, corresponding to the eigenvalues $\lambda_1$ and $\lambda_2$ ($\lambda_1$ and $\lambda_2$ can be conjugate), then
\[
M(t)=(q_1e^{\lambda_1 t},q_2e^{\lambda_2 t})
\]
is a fundamental matrix of Eq. (\ref{Constant}).
\end{theorem}

{\bf Proof.} From above discussion, we know $q_1e^{\lambda_1 t},q_2e^{\lambda_2 t}$ are two solution of Eq. (\ref{Constant}). Thus, $M(t)$ is a solution matrix of Eq. (\ref{Constant}). Moreover, by using the linear independence of $q_1$ and $q_2$, we have
\[
{\rm ddet}M(0)={\rm ddet}(q_1,q_2) \neq 0.
\]
Therefore, $M(t)$ is a fundamental matrix of Eq. (\ref{Constant}).

Now we need a lemma from (\cite{Baker} Proposition 2.4).

\begin{lemma}\label{Lem6.1}
 Suppose that $\lambda_1,\lambda_2,\cdots,\lambda_r$ are distinct eigenvalues for $A$, no two of which are
similar, and let $q_1,\cdots,q_r$ be corresponding eigenvectors. Then $q_1,\cdots,q_r$ are linearly independent.
\end{lemma}

Then by Lemma \ref{Lem6.1}, we have the corollary

\begin{corollary}
 If the matrix $A$ has two distinct eigenvalues $\lambda_1$ and $\lambda_2$ ($\lambda_1$ and $\lambda_2$ can not be conjugate), then
\[
M(t)=(q_1e^{\lambda_1 t},q_2e^{\lambda_2 t})
\]
is a fundamental matrix of Eq. (\ref{Constant}).
\end{corollary}

Now we need some results from (\cite{Bren} Theorem 11, Theorem 12).

\begin{lemma}\label{Lem6.2}
If $A$ is in triangular form, then every diagonal element is a
characteristic root.
\end{lemma}

\begin{lemma}\label{Lem6.3}
Let a matrix of quaternion be in triangular form. Then the only
characteristic roots are the diagonal elements (and the numbers similar to them).
\end{lemma}

\begin{lemma}\label{Lem6.4}
  Similar matrices have the same characteristic roots.
  \end{lemma}

\begin{lemma}\label{Lem6.5}
(\cite{Bren} Theorem 2)   Every matrix of quaternion can be transformed into triangular form by a unitary matrix.
\end{lemma}

\noindent {\bf Example 3} Find a fundamental matrix of the following QDEs
\begin{equation}
\dot {x}= \left ( \begin{array}{ll}
 \mathbf{i}  &  \mathbf{j}
 \\
 0  &  \mathbf{i+j}
\end{array}
\right )x,\,\,\,\, x=(x_1,x_2)^T.
\label{Ex6.2}
\end{equation}

{\bf Answer:} From Lemma \ref{Lem6.3} and Lemma \ref{Lem6.4}, we see that $\lambda_1=\mathbf{i}$ and $\lambda_2=\mathbf{i}+\mathbf{j}$.
To find the eigenvector of $\lambda_1=i$, we consider the following equation
\[
Aq=q\lambda_1,
\]
that is
\begin{equation}
 \left\{ \begin{array}{ccc}
\mathbf{ i}  q_1 + \mathbf{j}q_2 & =& q_1 \mathbf{i}
 \\
(\mathbf{i}+\mathbf{j}) q_2 &=& q_2 \mathbf{i}
\end{array}
\right.
\label{EV1}
\end{equation}
From the second equation of (\ref{EV1}), we can take $q_2=0$. Substituting it into the first equation of (\ref{EV1}), we obtain
$q_1=1$. So we can take one eigenvector as
\[
\nu_1=\left ( \begin{array}{ll}
 q_1
 \\
q_2
\end{array}
\right )
=\left ( \begin{array}{ll}
1
 \\
0
\end{array}
\right )
\]
To find the eigenvector of $\lambda_2=\mathbf{i}+\mathbf{j}$, we consider the following equation
\[
Aq=q\lambda_2,
\]
that is
\begin{equation}
 \left\{ \begin{array}{ccc}
 \mathbf{i}  q_1 + \mathbf{j} q_2 & =& q_1 (\mathbf{i}+\mathbf{j})
 \\
(\mathbf{i}+\mathbf{j}) q_2 &=& q_2 (\mathbf{i} + \mathbf{j})
\end{array}
\right.
\label{EV2}
\end{equation}
We can take one eigenvector as
\[
\nu_2=\left ( \begin{array}{ll}
 q_1
 \\
q_2
\end{array}
\right )
=\left ( \begin{array}{ll}
1
 \\
1
\end{array}
\right )
\]
Since
\[
\begin{array}{lll}
{\rm ddet}(\nu_1,\nu_2)
&=&
{\rm ddet}\left ( \begin{array}{ll}
 1 & 1
 \\
0 &  1
\end{array}
\right )
=
\displaystyle {\rm rdet}\Big[\left ( \begin{array}{ll}
 1 & 1
 \\
0 &  1
\end{array}
\right )\left ( \begin{array}{ll}
 1 & 0
 \\
1 & 1
\end{array}
\right )
\Big]
\neq0,
\end{array}
\]
 the eigenvectors $\nu_1$ and $\nu_2$ is linearly independent. Taking
 \[
 M(t)=(\nu_1 e^{\lambda_1 t},\nu_2 e^{\lambda_2 t})=\left ( \begin{array}{ll}
 e^{\mathbf{i} t}  & e^{(\mathbf{i}+\mathbf{j}) t}
 \\
0  &  e^{(\mathbf{i}+\mathbf{j}) t}
\end{array}
\right ),
 \]
 from Theorem \ref{Th6.1}, $M(t)$ is a fundamental matrix. In fact, by the definition of fundamental matrix, we can verify that $M(t)$ is a fundamental matrix of Eq. (\ref{Ex6.2}), which is consistent with Theorem \ref{Th6.1}. Now we verify the fundamental matrix as follows. First, we show that $M(t)$ is a solution matrix of Eq. (\ref{Ex6.2}).
 Let $\phi_1(t)=\nu_1 e^{\lambda_1 t}$ and $\phi_2(t)=\nu_2 e^{\lambda_2 t}$, then
  \[
\dot{\phi}_1(t)
=\left ( \begin{array}{ll}
1
 \\
0
\end{array}
\right )\mathbf{i}  e^{\mathbf{i} t}
=
\left ( \begin{array}{ll}
 \mathbf{i}  &   \mathbf{j}
 \\
 0  &  \mathbf{i}+\mathbf{j}
\end{array}
\right )
\left ( \begin{array}{ll}
 1
 \\
0
\end{array}
\right )  e^{\mathbf{i} t}
=
\left ( \begin{array}{ll}
 \mathbf{i}  &   \mathbf{j}
 \\
 0  &  \mathbf{i}+\mathbf{j}
\end{array}
\right ) \phi_1(t),
\]
which implies that $\phi_1(t)$ is a solution of Eq. (\ref{Ex6.2}).
Similarly, let $\phi_2(t)=\left ( \begin{array}{ll}
e^{(\mathbf{i}+\mathbf{j}) t}
 \\
e^{(\mathbf{i}+\mathbf{j}) t}
\end{array}
\right )$, then
\[
\dot{\phi}_2(t)=\left ( \begin{array}{ll}
1
 \\
1
\end{array}
\right )(\mathbf{i}+\mathbf{j})e^{(\mathbf{i}+\mathbf{j}) t}
=
\left ( \begin{array}{ll}
 \mathbf{i}  &   \mathbf{j}
 \\
 0  &  \mathbf{i}+\mathbf{j}
\end{array}
\right )
\left ( \begin{array}{ll}
1
 \\
1
\end{array}
\right )e^{(\mathbf{i}+\mathbf{j}) t}
=
\left ( \begin{array}{ll}
 \mathbf{i}  &   \mathbf{j}
 \\
 0  &  \mathbf{i}+\mathbf{j}
\end{array}
\right ) \phi_2(t),
\]
 which implies that $\phi_2(t)$ is another solution of Eq. (\ref{Ex6.2}). Therefore,  $M(t)=(\phi_1(t),\phi_2(t))^T$ is a solution matrix of Eq. (\ref{Ex6.2}).

Secondly, by Theorem \ref{Th5.2}, and taking $t_0=0$,
\[
\rm{ddet} M(t_0)=\rm{ddet} M(0) ={\rm ddet}(\nu_1,\nu_2) \neq0,
\]
then $\rm{ddet} M(t)\neq0$. Therefore, $M(t)$ is a fundamental matrix of Eq. (\ref{Ex6.2}).

\section{Conclusion and Discussion}

\noindent {\bf Conclusion}

This paper established a systematic basic theory for two-dimensional linear QDEs, which have many applications in the real world such as quantum mechanics, differential geometry, Kalman
filter design, attitude dynamics, fluid mechanics, and so on. For the sake of  understanding, we focus on the 2-dimensional system. However, in fact, our basic results can be easily extended to any arbitrary $n$-dimensional QDEs.
Due to the non-commutativity of quaternion algebra, there are four profound differences between QDEs and ODEs.

1. Due to the non-commutativity of the quaternion algebra, the algebraic structure of the solutions to QDEs is not a linear vector space. It is actually a right-free module. 

2. {\em Wronskian} of ODEs is defined by Caley determinant. However, since Caley determinant depends on the expansion of $i$-th row and $j$-th column of quaternion, different expansions can lead to different results. {\em Wronskian} of ODEs can not be extended to QDEs. It is necessary to define the {\em Wronskian} of QDEs by a novel method (see Section 4). We use double determinants $\rm {ddet}M=\displaystyle \rm{rddet}(MM^{+})$ to define {\em Wronskian} of QDEs. We remark that the {\em Wronskian} of QDEs can also be defined by $\rm {ddet}M=\displaystyle \rm{rddet}(M^{+}M)$.

3. Liouville formula for QDEs and ODEs are different.

4. For QDEs, it is necessary to treat the eigenvalue problems with left and right, separately. This is a large-difference between QDEs and ODEs.

\section{Conflict of Interests}

The authors declare that there is no conflict of interests
regarding the publication of this article.

\section{Acknowledgement}
   The authors would express their great gratitude to Professor Weinian Zhang for his valuable comments and discussion when they prepared this paper.

\end{document}